%% file: multiscale.tex
\documentclass[final,1p,12pt]{elsarticle}

\usepackage{amsfonts,amsopn,amsmath,amssymb,amsthm}
\usepackage{mathrsfs}
\usepackage{color}
\usepackage{booktabs}
\usepackage{nicefrac}
\usepackage{multirow}
\usepackage{enumerate}
\usepackage{tikz}
\usetikzlibrary{shapes.geometric, patterns, spy}
\usepackage{pgfplots}
\usepackage{pgfplotstable}
\usepgfplotslibrary{groupplots}
\usepackage{graphicx}

\usepackage{placeins}
\pgfplotsset{compat=1.18}
\usepackage{subcaption}
\usepackage{adjustbox}

\usepackage{mathtools}
\usepackage{pdfpages}

\usepackage[hidelinks]{hyperref}
\usepackage{xcolor}
\hypersetup{
    colorlinks,
    linkcolor={red!50!black},
    citecolor={blue!50!black},
    urlcolor={blue!80!black}
}
\usepackage{cleveref}
\usepackage{url}
\graphicspath{{Fig/}}

\journal{arxiv}

\newcommand{\NN}{\mathbb{N}}
\newcommand{\RR}{\mathbb{R}}
\newcommand{\CC}{\mathbb{C}}
\newcommand{\BB}{\mathbb{B}}

\renewcommand{\O}{\mathcal{O}}
\newcommand{\V}{\mathcal{V}}
\newcommand{\D}{\mathcal{D}}
\newcommand{\Jx}{J_x}
\newcommand{\Jy}{J_y}
\newcommand{\Dincl}[1]{\mathcal{D}_{#1 \times #1}}

\newcommand{\M}{\mathcal{M}}

\newcommand{\PpBd}{\mathcal{P}_p}
\newcommand{\Oj}{\Omega_j}
\newcommand{\T}{\mathcal{T}}
\newcommand{\Th}{\mathcal{T}_h}

\newcommand{\VFEM}{V^{h,p}}
\newcommand{\VFEMO}{V_{0}^{h,p}}

\newcommand{\VACMSh}{V_{A}^{h,p}(\D)}
\newcommand{\VACMS}{V_{A}(\D)}

\newcommand{\SFEM}{\mathbb{S}_{F}}
\newcommand{\NFEM}{N_{F}}
\newcommand{\NFEMj}{N_{j}}
\newcommand{\NFEMe}{n_{e}}

\newcommand{\uFEM}{u^{h,p}}
\newcommand{\uGround}{u^{h,10}}
\newcommand{\vFEM}{v^{h,p}}
\newcommand{\uACMSh}{u^{h,p}_A}
\newcommand{\vACMSh}{v^{h,p}_A}
\newcommand{\uACMS}{u_A}
\newcommand{\vACMS}{v_A}

\newcommand{\SACMS}{\mathbb{S}_{A}}
\newcommand{\NACMS}{N_{A}}

\newcommand{\Cext}{C_{ext}}
\newcommand{\Ceig}{C_{eig}(\Ie)}
\newcommand{\Cstab}{C_{\kappa}}

\newcommand{\Ie}{I_{\E}}

\newcommand{\uvecFEM}{\mathbf{u}_{F}}

\newcommand{\gvecFEM}{\mathbf{g}_{F}}
\newcommand{\uvecACMS}{\mathbf{u}_A}

\newcommand{\gvecACMS}{\mathbf{g}_A}
\newcommand{\xx}{\mathbf{x}}

\newcommand{\eigh}{\lambda_{e,i}^h}
\newcommand{\tauh}{\tau^{h,p}}
\newcommand{\vh}{v^{h,p}}
\newcommand{\tauhe}{\tau_{e,i}^{h,p}}
\newcommand{\phiPh}{\varphi_h^q}
\newcommand{\etah}{\eta^{h,p}}

\newcommand{\sesH}{\mathcal{C}} 

\newcommand{\taue}{\tau_i^e}
\newcommand{\eig}{\lambda^e_i}
\newcommand{\phiP}{\varphi^q}

\newcommand{\E}{\mathcal{E}}
\newcommand{\Ve}{\mathcal{V}}

\newcommand{\ExtG}{E^\Gamma} 

\newcommand{\ExtLoc}[1]{E^{#1}}

\newcommand{\ExtGh}{E_h^\Gamma} 

\newcommand{\ExtLoch}[1]{E_h^{#1}}

\newcommand{\HOj}{H^1_0(\Oj)}
\newcommand{\Hj}{H^1(\Oj)}
\newcommand{\HO}{H^1(\Omega)}
\newcommand{\LO}{L^2(\Omega)}

\newcommand{\HtrG}{H^{1/2}(\Gamma)}
\newcommand{\HtrBOj}{H^{1/2}(\partial\Oj)}
\newcommand{\HtrBO}{H^{1/2}(\partial\Omega)}
\newcommand{\HOe}{H^1_0(e)}

\newcommand{\cplu}{\iota}

\newtheorem{theorem}{Theorem}

\theoremstyle{plain}

\theoremstyle{definition}

\newtheorem{remark}[theorem]{Remark}

\pgfplotscreateplotcyclelist{philipcolors}{%
gray,dashed,very thick,mark size=3pt\\%
gray,dotted,very thick,mark size=3pt\\%
brown,dashed,very thick,mark size=3pt\\%
brown,dotted,very thick,mark size=3pt\\%
violet,very thick,mark size=3pt\\%
teal,very thick,mark size=3pt\\%
orange,very thick,mark size=3pt\\%
cyan,very thick,mark size=3pt
\\
violet,very thick,mark size=3pt, mark=*, mark options={scale=0.5}\\%
teal,very thick,mark size=3pt, mark=*, mark options={scale=0.5}\\%
orange,very thick,mark size=3pt, mark=*, mark options={scale=0.5}\\%
cyan,very thick,mark size=3pt, mark=*, mark options={scale=0.5}\\
purple,very thick,mark size=3pt, mark=*, mark options={scale=0.5}\\
pink,very thick,mark size=3pt, mark=*, mark options={scale=0.5}\\
magenta,very thick,mark size=3pt, mark=*, mark options={scale=0.5}
\\
violet,very thick,mark size=3pt, mark=*, mark options={scale=0.3}\\%
teal,very thick,mark size=3pt, mark=*, mark options={scale=0.3}\\%
orange,very thick,mark size=3pt, mark=*, mark options={scale=0.3}\\%
cyan,very thick,mark size=3pt, mark=*, mark options={scale=0.3}\\
purple,very thick,mark size=3pt, mark=*, mark options={scale=0.3}\\
pink,very thick,mark size=3pt, mark=*, mark options={scale=0.3}\\
magenta,very thick,mark size=3pt, mark=*, mark options={scale=0.3}\\
}

\begin{document}

\begin{frontmatter}



\title{High-order discretized ACMS method for the simulation of finite-size two-dimensional photonic crystals}

\affiliation[label1]{organization={Department of Applied Mathematics, University of Twente},
            addressline={P.O.~Box~217}, 
            city={Enschede},
            postcode={7500 AE}, 
            country={The Netherlands}}
\affiliation[label2]{organization={Delft Institute of Applied Mathematics, Delft University of Technology},
             addressline={Mekelweg~4},
             city={Delft},
             postcode={2628 CD},
             country={The Netherlands}}
\affiliation[label3]{organization={Department of Mathematics, University of Hamburg},
             addressline={Bundesstraße~55},
             city={Hamburg},
             postcode={20146},
             country={Germany}}

\author[label1]{Elena Giammatteo\fnref{fn1}} 
\ead{e.giammatteo@utwente.nl}
\author[label2]{Alexander Heinlein} 
\ead{a.heinlein@tudelft.nl}
\author[label3]{Philip L. Lederer} 
\ead{philip.lederer@uni-hamburg.de}
\author[label1]{Matthias Schlottbom\corref{cor1}} 
\ead{m.schlottbom@utwente.nl}
\cortext[cor1]{Corresponding author}
\fntext[fn1]{The research of this author was supported by the Dutch Research council (NWO) via grant OCENW.GROOT.2019.071.}

\begin{abstract}
The computational complexity and efficiency of the approximate mode component synthesis (ACMS) method is investigated for the two-dimensional heterogeneous Helmholtz equations, aiming at the simulation of large but finite-size photonic crystals. The ACMS method is a Galerkin method that relies on a non-overlapping domain decomposition and special basis functions defined based on the domain decomposition. 
While, in previous works, the ACMS method was realized using first-order finite elements, we use an underlying $hp$--finite element method.
We study the accuracy of the ACMS method for different wavenumbers, domain decompositions, and discretization parameters. Moreover, the computational complexity of the method is investigated theoretically and compared with computing times for an implementation based on the open source software package NGSolve.
The numerical results indicate that, for relevant wavenumber regimes, the size of the resulting linear systems for the ACMS method remains moderate, such that sparse direct solvers are a reasonable choice. Moreover, the ACMS method exhibits only a weak dependence on the selected domain decomposition, allowing for greater flexibility in its choice. Additionally, the numerical results show that the error of the ACMS method achieves the predicted convergence rate for increasing wavenumbers. Finally, to display the versatility of the implementation, the results of simulations of large but finite-size photonic crystals with defects are presented.
\end{abstract}




\end{frontmatter}



\section{Introduction} \label{Sec:Introduction}
Photonic crystals are nanostructures consisting of periodically repeating building blocks, the so-called unit cells, which usually consist of different materials, modeled by discontinuous material parameters \cite{AdhikaryUppuVos2020,Joannopolous2008}.
Photonic crystals find many applications, such as enhanced absorption in solar cells and thin films \cite{CallahanHorowitzAtwater2013,ZhouBiswas2008,SharmaVegtVos2021}, or the control of spontaneous emission of quantum emitters \cite{KoenderinkBechgerVos2002,LeistikowMoskVos2011,LodahlDrielVos2004,OgawaImadaYoshimoto2004}.
A mathematical treatment of photonic crystals can be found in \cite{DorflerLechleiterPlumSchneiderWieners2011,Kuchment2016}.

The governing equations for light propagation in photonic crystals are Maxwell's equations, which may be reduced, under certain conditions, to the heterogeneous Helmholtz equations for 
transverse-magnetic and transverse-electric polarizations \cite{Joannopolous2008}.
In view of the complex behavior of electromagnetic waves in photonic crystals, numerical approaches for the simulation of their optical properties are required, which we briefly discuss next. 
%

Finite difference methods \cite{BuschFreymannLindenMingaleevTkeshelashviliWegener2007,HoIvanovEnginEtAl2011,MinkovSavona2015,TafloveHagness2005,Yee1966} are primarily used for structured grids, which allow for computationally highly efficient implementations. Using the regular sparsity structure of the resulting linear systems, problems up to tens of millions to one hundred million degrees of freedom can still be solved using sparse direct solvers~\cite{AmestoyBrossierButtari2016,OpertoAmestoyButtari2023}. 
Second order finite difference methods may suffer from a so-called pollution effect, i.e., large dispersion errors, which require much finer meshes than expected for a given wavenumber. 
Higher order finite difference methods weaken such strong assumption on the mesh, see, e.g., \cite{TongFangZhao2023} for a discussion.
However, the required smoothness of the fields may not be satisfied in the aforementioned application, e.g., due to discontinuous material parameters. 
Moreover, stair casing may occur for curved geometries, which requires non-trivial modifications for a robust and accurate implementation \cite{Nicolaides2004}.

Finite element methods (FEM) are particularly well-suited for handling discontinuous material parameters due to their high geometric flexibility, but they also suffer from a pollution effect, which, in the FEM context, shows up in the form of wavenumber-dependent constants in the error analysis \cite{BabuskaIhlenburg1995,BabuskaSauter1997}. 
As a result, linear finite elements generally require excessive mesh resolution for growing wavenumbers, whereas high-order methods allow for overcoming the pollution effect; see, e.g., \cite{ChaumontFreletNicaise2020,ChenLuXu2013,MelenkSauter2011,MelenkParsaniaSauter2013,LafontaineSpence2022}. 
Let us also refer to interior penalty (high-order) discontinuous Galerkin methods \cite{FengWu2009,FengWu2011} for the homogeneous Helmholtz equation.
Additionally, boundary integral equation methods have been developed in \cite{StrauszerFariaPerez2023} to solve plane wave scattering problems by periodic arrays of two-dimensional obstacles.

For large nanostructures, finite element methods lead to large linear systems, because of the required partitioning of the unit cells into sufficiently small elements. Therefore, solving these systems using a direct solver may become infeasible.
Additionally, solving the resulting systems using iterative methods is also challenging for large wavenumbers:~cf., e.g.,~\cite{GanderZhang2019,GanderGrahamSpence2015}.
To reduce the size of the linear systems, approaches that employ problem-adapted basis functions instead of piecewise polynomials have been employed.
Trefftz discontinuous Galerkin methods~\cite{GittelsonHiptmairPerugia2009,HiptmairMoiolaPerugia2016,YuanGong2023} use solutions of the differential equation (plane waves) for the local approximation of the fields. For piecewise constant materials, these are often known, but for non-constant materials, the computation of such functions becomes more involved. In~\cite{LehrenfeldStocker2023}, polynomial discrete plane wave functions have been constructed numerically from a piecewise polynomial discretization. From all the degrees of freedom per element when considering a standard discontinuous Galerkin method, some are removed, while others are agglomerated into (polynomial) Trefftz ansatz functions. A detailed comparison for the latter is also given in \cite{LedererLehrenfeldStocker2024}.

For infinite periodically repeating nanostructures, the solutions are given in terms of Bloch modes \cite{DorflerLechleiterPlumSchneiderWieners2011, LuVegtZu2017}. In \cite{BrandsmeierSchmidtSchwab2011}, a conforming finite element method that uses Bloch modes has been developed for 2D  photonic crystals that are infinite in one direction and surrounded by an external medium in the other. The conformity requirement is satisfied by modulating the Bloch modes with polynomial functions, which makes the assembly of the corresponding matrices involved, but independent of the size of the structure.
In~\cite[Chapter~5]{Kozon2023}, photonic crystals which are finite in two-dimensions, are discretized using a discontinuous Galerkin method that treats the crystal as one element and that uses Bloch modes as basis functions within the crystal. Dropping the conformity requirement and exploiting properties of the Bloch modes greatly simplifies the assembly, yielding only a weak dependence on the number of unit cells; the dependence is linear instead of quadratic for crystals consisting of the same amount of unit cells in each direction. 

When the periodic arrangement of the unit cells is broken, the aforementioned methods \cite{BrandsmeierSchmidtSchwab2011,Kozon2023,StrauszerFariaPerez2023} are, however, not directly applicable. Such cases occur, e.g., in waveguides \cite{Joannopolous2008} or naturally in truly manufactured nanostructures due to unavoidable manufacturing flaws. For direct numerical simulations of real nanostructures and comparisons to the design structure using a discontinuous Galerkin method, see \cite{CorbijnSchlottbomVegtVos2024}.

Multiscale methods \cite{ChenHouWang2023,GiammatteoHeinleinSchlottbom2024,MaAlberScheichl2023,Peterseim2017,PeterseimVerfurth2020} employ (local) problem adapted basis functions. As such they aim to provide accurate approximation with only a few degrees of freedom, and at the same time, they can handle non-constant material properties and complex geometries. In practice, however, such methods usually also rely on an underlying discretization for their numerical realization.

This work focuses on the computational aspects and implementation of the approximate mode component synthesis (ACMS) method to solve two-dimensional heterogeneous Helmholtz equations. We refer to~\cite{GiammatteoHeinleinSchlottbom2024} for a theoretical analysis and preliminary numerical studies. The method has originally been introduced in~\cite{Hetmaniuk2010} and further theoretically analyzed in~\cite{Hetmaniuk2014} for heterogeneous diffusion problems.
The ACMS method relies on a combination of a non-overlapping domain decomposition approach, which allows for straightforward parallelization over the number of unit cells during assembly, and the construction of problem-adapted basis functions.
We mention that the ACMS method does not rely on (local) periodicity assumptions, allowing it to easily handle waveguides, defects, and other perturbations in the crystal structure. This flexibility broadens the method's applicability, including its use in the numerical simulation of quasicrystals; cf.~\cite{VardenyNahataAgrawal2013}.
The computational experiments in~\cite{GiammatteoHeinleinSchlottbom2024} considered relatively small examples and employed linear finite elements, which, as mentioned above, are susceptible to the pollution effect. Consequently, also the discrete realization of the ACMS method can potentially suffer from the pollution effect. Moreover, the efficient implementation of the ACMS method and its computational complexity has not been investigated in \cite{GiammatteoHeinleinSchlottbom2024}; however, an efficient parallel implementation of the ACMS method employing linear finite elements for two-dimensional heterogeneous diffusion problems has been presented in~\cite{Heinlein2015}.

Before going into a more detailed description of the ACMS method, we outline our main contributions, aimed at addressing the research gaps mentioned above. 
Our contributions include a detailed computational complexity analysis and an efficient implementation of the ACMS method using an underlying $hp$--finite element method.
The efficiency of our implementation, i.e., assembly of linear systems and their solution using sparse direct solvers, is due to structural similarities to usual $hp$--finite element methods. In particular, the usual local degree $p$ polynomials are replaced by problem-adapted basis functions; and when fixing the domain decomposition, the ACMS method allows for increasing the number of local basis functions for each ACMS element (subdomain), similar to increasing the polynomial degree in $hp$--finite elements. 
The rather small number of basis functions required to achieve engineering relevant tolerances in practical wavenumber regimes yields linear systems of moderate size, allowing the use of sparse direct solvers.
Additionally, since our implementation is based on underlying high-order finite elements, we can numerically investigate the wavenumber dependence of the method. This is particularly important for the Helmholtz equations, because Galerkin methods like the ACMS, require sufficient resolution in order to be well defined; see, e.g., \cite{GrahamSauter2019}. 
Moreover, we study the influence on the computational complexity as well as on the accuracy of the choice of the domain decomposition and the number of modes per subdomain.
Finally, we show that, by suitable choices of the discretization parameters, we can simulate large two-dimensional nanostructures.
Our implementation employs the open source package NGSolve~\cite{NGSolve} and the code is publicly available at~\cite{NGSolveACMS}.

The remainder of the paper is structured as follows. We first introduce the Helmholtz equation and the essential notation in \Cref{Sec:ProblemFormulation}. Then, in \Cref{Sec:ACMS}, we give a summary of the continuous ACMS method and recall relevant analytical results. In \Cref{Sec:DiscreteACMS}, we describe the numerical realization and computational cost of the ACMS method when employing an underlying $hp$--FEM discretization. Finally, in \Cref{Sec:NumericalResults}, we present numerical results to highlight the method's flexibility in the choice of the domain decomposition and the benefit of choosing a high-order approximation of the ACMS basis functions. We conclude with final remarks in~\Cref{sec:discussion}.

\section{Notation and model equation} \label{Sec:ProblemFormulation}

In the coming section, we present the basic notation along with the Helmholtz equation.

\paragraph{Notation}
Given a connected polygonal domain $\Omega \subset \RR^2$  with piecewise $C^2$ boundary and strictly convex angles, we denote with $\LO$ the Lebesgue space of square-integrable functions $u, v : \Omega \rightarrow \CC$ with inner product 
\begin{align*}
(u, v)_{\Omega} = \int_{\Omega} u \overline{v} \ d\xx   ,
\end{align*}
and with $H^1(\Omega)$ the usual Sobolev space of functions in $\LO$ with square-integrable weak derivatives. Additionally, let $\HtrBO$ be the trace space of functions in $\HO$ and let $L^{\infty}(\Omega)$ be the space of essentially bounded measurable functions. 
Corresponding notation is also used for other measurable sets besides $\Omega$. We indicate with $\| \cdot \|_V$ the corresponding norm on a function space $V$.

\paragraph{Problem formulation}
We consider the heterogeneous Helmholtz problem:
\begin{gather}
\begin{aligned}\label{eq:model_strong}
-{\rm div}(a \nabla u) - \kappa^2 u &= f \quad \text{in } \Omega , \\
 a \partial_n u - \cplu  \omega \beta u &=g  \quad  \text{on } \partial \Omega  . 
\end{aligned}
\end{gather}
The coefficient functions $a,c \in L^{\infty}(\Omega)$ describe material properties of the background medium occupying $\Omega$ and are such that $a \geq a_0 >0$ and $c \geq c_0 >0$ for some constants $ a_0, c_0 \in \RR$. Together with the positive angular frequency denoted by $\omega$, we then obtain the wavenumber $\kappa=\omega /c$.
The real-valued function $\beta \in L^{\infty} (\partial \Omega)$  is related to transmission and reflection of the unknown function $u$ on the boundary $\partial \Omega$, and we assume that either $\beta \geq \beta_0 >0 $ or $\beta \leq \beta_0 <0 $, for a given $\beta_0 \in \RR$.
Impedance boundary conditions are described by the function $g \in \HtrBO$, and any interior sources may be modeled by the function $f \in \LO$.
In the current work we focus on the case of a homogeneous interior source, i.e., $f \equiv 0$, and we refer to \cite{GiammatteoHeinleinSchlottbom2024} for a discussion on the method in the general case of a non-zero source term.

To obtain a weak formulation of the Helmholtz problem, let us introduce the sesquilinear form $\sesH: \HO \times \HO \to \CC$ defined by
\begin{align*} 
\sesH(u,v) = (a \nabla u, \nabla v)_{\Omega} - (\kappa^2 u,v)_{\Omega} - \cplu (\omega \beta u,v)_{\partial \Omega}  ,
\end{align*}
and the antilinear functional $G:\HO \to \CC$ defined by
\begin{align*} 
    G(v) = (g,v)_{\partial \Omega}.
\end{align*}
The weak form of the Helmholtz problem~\eqref{eq:model_strong} then reads
\begin{align} \label{eq:HelmholtzVariational}
    &\text{Find } u \in \HO: \quad \sesH(u,v) = G(v), \quad \text{for all } v \in \HO.
\end{align}
With the aforementioned regularity assumptions on the coefficients and the functions, the Helmholtz problem is well-posed and the weak formulation \eqref{eq:HelmholtzVariational} admits a unique solution; see~\cite[Theorem~2.4]{GrahamSauter2019}.

\section{The continuous ACMS method}\label{Sec:ACMS}
In the following section, we recap the ACMS method and recall some theoretical convergence results obtained in \cite{GiammatteoHeinleinSchlottbom2024}.
\paragraph{Domain Decomposition} 
We introduce by $\D = \{ \Oj, \ j=1,...,J \}$ a conforming decomposition of $\Omega$ into $J$ non-overlapping subdomains with piecewise smooth boundaries, and we define the interface $\Gamma$ of the domain decomposition as
\begin{align*} 
    \Gamma = \bigcup_{j=1}^J \partial\Oj  ,
\end{align*}
which also includes the boundaries of the domain $\Omega$.
Moreover, we introduce the set $\E$ of edges and the set $\V$ of vertices of the domain decomposition:
\begin{align*}
    \E &= \{e \subset \Gamma: \ e = \partial \Omega \cap \partial \Oj \ \text{ or } \ e = \partial \Oj \cap \partial \Omega_i , \text{ for some } i,j = 1,...,J \} , 
    \\
    \V &= \{q \in \Gamma: \ \exists e \in \E:  \ q \in \partial e \} .
\end{align*}

\paragraph{ACMS Spaces}
The proper formulation of the ACMS method relies on the solvability of local Helmholtz problems with Dirichlet boundary conditions \cite{GiammatteoHeinleinSchlottbom2024}. To that end, fix $j\in\{1,\ldots,J\}$, and let us consider the following eigenvalue problem: Find $b\in H^1_0(\Oj)\setminus\{0\}$ and $\lambda \in \RR$ such that
\begin{align*}
(a \nabla b, \nabla v)_{\Oj} =\lambda (\kappa^2 b,v)_{\Oj} \quad\text{for all } v \in  H^1_0(\Oj).
\end{align*}
The local solvability condition can then be conveniently characterized by requiring that $\lambda \neq 1$, for all $j\in\{1,\ldots,J\}$, which we assume in the following.
Under this assumption, and recalling that we assumed $f=0$, the solution $u$ to \eqref{eq:HelmholtzVariational} is determined by its trace on $\Gamma$.
The ACMS method now relies on the approximation of $u|_\Gamma$ using functions that are locally supported on $\Gamma$.
We will introduce the local functions associated to vertices and edges as follows. 
For all $q \in \V$, let $\phiP: \Gamma\to \RR$ be edgewise harmonic functions defined on $\Gamma$, i.e.,  $\phiP(q)=\delta_{q,r}$ for $q,r \in \V$, where $\delta$ is the Kronecker delta, and
\begin{align*}
\int_e \partial_e \phiP \partial_e \eta \ ds = 0 \quad \text{for all } \eta \in \HOe,
\end{align*}
where $\partial_e$ is the tangential derivative along the edge $e$. 
The vertex-based space is then defined as the linear combination of the functions $\phiP$:
\begin{align*}
V_{\Ve} = {\rm span}\{ \phiP: \ q \in \Ve\}.
\end{align*}
Additionally, let us define the edge modes as solutions to the following weak formulation of the edge-Laplace eigenvalue problems: for each $e\in\E$, find $(\taue,\eig)\in \HOe\times\RR$, $i \in \NN$, such that
\begin{align}\label{eq:def_edge_mode}
(\partial_e \taue,\partial_e \eta)_e =\eig ( \taue, \eta)_e, \quad \text{for all } \eta\in \HOe.
\end{align}
According to \cite[p.~415]{CourantHilbert1953}, the eigenvalues depend quadratically on the index of the mode, that is
\begin{align} \label{eq:eig_behaviour}
    \lambda_i^e \sim \Big( \frac{i \pi}{|e|} \Big)^2.
\end{align}
With a slight abuse of notation, we also denote by $\taue$ its extension by zero to the whole interface $\Gamma$. 
The edge-based space is then defined as the linear combination of the eigenmodes $\taue$ for all edges
\begin{align*}
V_{\E} =  {\rm span}\{ \taue: \ e \in \E, \ i \in \NN \} .
\end{align*}
It follows from \cite[Lemma~3.6]{GiammatteoHeinleinSchlottbom2024}, that $V_{\Ve} + V_{\E}$ is dense in the trace space
\begin{align*}
    \HtrG = \{v:\Gamma\to\CC: \  \forall j = 1,...,J, \ \exists u_j \in \Hj \mbox{ s.t. } u_{j}|_{\partial\Oj} = v|_{\partial\Oj}  \}.
\end{align*}

In order to approximate \eqref{eq:HelmholtzVariational}, we next extend functions in $\HtrG$ to functions in $\HO$. To that end, we first introduce the extensions $\ExtLoc{j}$ associated to the subdomains $\Oj$, and then the extension $\ExtG$ associated to the full interface $\Gamma$. 
For a given $\Oj \in \D$ and $\tau \in \HtrBOj$, let us define the local Helmholtz-harmonic extension operator $\ExtLoc{j}:\HtrBOj\to \Hj$ such that 
\begin{gather}
\begin{aligned}\label{eq:localHarmonicExtension}
    (a \nabla \ExtLoc{j} \tau, \nabla \eta)_{\Oj} - (\kappa^2\ExtLoc{j} \tau, \eta)_{\Oj} &= 0  , \quad\text{for all } \eta \in  \HOj  , \\
      (\ExtLoc{j} \tau )|_{\partial \Oj} &= \tau .
\end{aligned}
\end{gather}
The Helmholtz-harmonic extension operator is well-defined, which follows from the above assumption on the solvability of local Helmholtz problems, and bounded \cite{GiammatteoHeinleinSchlottbom2024}. 
Let us also introduce $\ExtG:\HtrG\to \HO$, which first restricts a given function $\tau \in \HtrG$ to the boundaries of the subdomains $\Oj$, for all $j = 1,\ldots,J$, and then applies the local Helmholtz-harmonic extension operator, i.e., $(\ExtG \tau)|_{\Oj} = \ExtLoc{j} (\tau|_{\partial \Oj})$. 
We note that a vertex basis function $\phiP$ is supported on all edges that share the vertex $q$ while its extension $\ExtG\phiP$ is supported on all subdomains $\Oj$ that share the vertex $q$. In contrast, the extension of an edge basis function $\ExtG \taue $ is supported on the two neighboring subdomains that share the edge $e$. 
Moreover we can write the solution $u$ of~\eqref{eq:HelmholtzVariational} in terms of its trace $u|_{\Gamma} \in \HtrG$ that is extended to $\HO$ with the operator that was just introduced: $u= \ExtG (u|_{\Gamma})$. 

Using the above extensions, we can now approximate the solution space using a combination of local spaces associated to vertex functions and edge modes:
\begin{align*} 
    V(\D) := \ExtG V_{\Ve} + \ExtG V_{\E} \subseteq \HO.
\end{align*}

We are now able to define the finite-dimensional ACMS space $\VACMS \subset V(\D)$.
For every edge, we select the eigenmodes associated to the $\Ie \in \NN$ smallest eigenvalues, assuming that they are non-decreasingly ordered, thus obtaining
\begin{align} \label{eq:approximation_space_finite}
    \VACMS := \ExtG V_{\Ve} + \ExtG \V_{\E}^{\Ie} , \quad \mbox{ with } \V_{\E}^{\Ie} = {\rm span}\{ \taue : \ e \in \E,  \ 1 \leq i \leq \Ie \} .
\end{align}

\paragraph{ACMS approximation of the Helmholtz equation} \label{Sec:ContinuousACMS}

The weak formulation \eqref{eq:HelmholtzVariational} of the Helmholtz equation can be approximated in the ACMS space: find $\uACMS \in \VACMS $ such that
\begin{align} \label{eq:ACMS_Weak_Formulation}
\sesH(\uACMS, \vACMS) = G(\vACMS), \quad \text{for all } \vACMS \in \VACMS   .
\end{align}
The well-posedness of the approximation problem \eqref{eq:ACMS_Weak_Formulation} relies on the smallness of the adjoint approximability constant $\sigma^*$. In \cite[Theorem~4.6]{GiammatteoHeinleinSchlottbom2024}, 
it has been shown that 
  $\sigma^* = \O(\Cstab \| \kappa\|_{\infty}^2/\Ie)$ if the coefficients $a$ and $\beta$ are Lipschitz continuous in a neighborhood of the interface $\Gamma$.
Here, $\Cstab$ is the stability constant of the Helmholtz problem, i.e., a bound on the solution operator associated with \eqref{eq:HelmholtzVariational}, mapping $f\in L^2(\Omega)$ to $u\in H^1(\Omega)$, for $g=0$.
We note that $\Cstab$ in general depends on $\kappa$, and we refer to \cite{GrahamSauter2019} for a discussion.
Hence, $\sigma^*$ can be made small by using sufficiently many edge modes.
Moreover, under these regularity assumptions on $a$ and $\beta$, it has been shown in \cite[Theorem~4.7]{GiammatteoHeinleinSchlottbom2024} that
\begin{align}  \label{eq:L2_error_bound_kappa_dep}
        \| u - \uACMS \|_{\LO} 
        &\leq C \frac{ \max_{e\in\E}|e|^3 \Cstab \| \kappa\|_{\infty}^2 }{\Ie^3} \sum_{e \in \E} \| u \|_{H^3(e)},
\end{align}
where $C$ is a constant independent of $\kappa$.
These results indicate that the wavenumber $\kappa$ influences the number of edge modes required to obtain a given $L^2$--error of the ACMS approximation against the exact solution. 
Since other factors enter the estimate, the relation between edge modes and wavenumber in this bounds will be further investigated numerically; see \Cref{test:numex1,test:numex2}.

\section{The discrete ACMS method}\label{Sec:DiscreteACMS}
The approximation space~\eqref{eq:approximation_space_finite} is spanned by extensions of the vertex basis functions and of the edge modes defined in \eqref{eq:def_edge_mode}. In practice, though, we need to numerically approximate the eigenproblems and the harmonic extension operators. We now discuss the numerical realization of the ACMS method when employing an $hp$--FEM discretization.
\paragraph{Triangulation}
For the sake of simplicity, we suppose that all edges are line segments, although the coming discussion would follow similarly in case of curved edges.
Let $\Th =\{ T_1,...,T_N\}$ be a triangulation of the domain $\Omega$ without hanging nodes, with mesh size
\begin{align*}
    h = \max_{T \in \T } h_T ,
\end{align*}
where $h_T$ is the diameter of $T$; see, e.g., \cite{BrennerScott2008}. If not explicitly mentioned differently, we use a quasi uniform triangulation. Additionally, we want the mesh to be conforming with the domain decomposition, namely, we want every triangle $T \in \Th$ to be contained exactly in one subdomain $T \subset \Oj$.

\paragraph{FEM approximation of the Helmholtz equation}
For the sake of the discussion, let us present a finite element approximation of the Helmholtz equation. \\
Let $\PpBd$ be the set of all polynomials in two variables of degree less or equal to $p$. Then, the usual finite element space $\VFEM(\Omega)$ defined on the triangulation of the domain $\Omega$ consisting of piecewise polynomials of degree $p$ is given by
\begin{align*}
    \VFEM(\Omega) = \{ v \in C^0(\overline{\Omega}): \ v|_{T} \in \PpBd, \ \forall \ T \in \Th \} ;
\end{align*}
see, e.g., \cite{BrennerScott2008}. 
The finite element approximation of the weak formulation~\eqref{eq:HelmholtzVariational} reads as:  find $\uFEM \in \VFEM(\Omega)$ such that
\begin{align} \label{eq:FEM_Weak_Formulation}
\sesH(\uFEM, \vFEM) = G(\vFEM), \quad \text{for all } \vFEM \in \VFEM(\Omega).
\end{align}
By choosing a basis for $\VFEM(\Omega)$, equation \eqref{eq:FEM_Weak_Formulation} reads 
\begin{align*} 
    \SFEM \uvecFEM = \gvecFEM  ,
\end{align*}
where $\uvecFEM \in \RR^{\NFEM}$ is the solution vector, $\SFEM \in \RR^{\NFEM \times \NFEM}$ is the system matrix and $\gvecFEM \in \RR^{\NFEM}$ is the right-hand side vector, and $\NFEM = \dim \VFEM(\Omega)$.

\subsection{Setup of the discrete ACMS method}
In order to define the discrete counterparts of the elements in $\VACMS$, we introduce auxiliary spaces on a generic set $\M \subseteq \overline{\Omega}$ in the following way:
\begin{align*}
    \VFEM(\M) &= \{ v|_{\M} : \ v \in \VFEM(\Omega) \}  , \\
    \VFEMO(\M) &= \{ v|_{\M} : \ v \in \VFEM(\Omega), \ v|_{\partial \M} = 0 \}  .
\end{align*}
The discrete vertex basis function of a generic vertex $q \in \V$ is referred to as $\phiPh \in \VFEM(\Gamma)$ and is such that $\phiPh(r) = \delta_{q,r}$ for all $q,r \in \V$. Moreover, the function is piecewise harmonic on the edges in the discrete sense, meaning that it satisfies
\begin{align} \label{eq:vertex_harmonic}
     ( \partial_e \phiPh , \partial_e \etah)_e  &= 0  , \quad\text{for all } \etah \in  \VFEMO(e), \ \text{for all} \ e \in \E   .
\end{align}
The discrete edge modes are computed by numerically solving the eigenvalue problem~\eqref{eq:def_edge_mode}: for $e\in\E$ and $i=1,...,\Ie$, find $(\tauhe,\eig)\in \VFEMO(e) \times \RR$ such that
\begin{align}\label{eq:edge_mode_discrete}
(\partial_e \tauhe, \partial_e \etah)_e =\eigh ( \tauhe, \etah)_e  , \quad\text{for all } \etah \in \VFEMO(e)  .
\end{align}
Again, with a slight abuse of notation we also denote by $\tauhe$ its extension by zero to the whole interface $\Gamma$.
As in~\eqref{eq:localHarmonicExtension}, we introduce the discrete Helmholtz-harmonic extension $\ExtLoch{j}: \VFEM(\partial \Oj) \to \VFEM(\Oj)$ such that:
\begin{gather}
\begin{aligned}\label{eq:localHarmonicExtensionDiscrete}
    (a \nabla \ExtLoch{j} \tauh, \nabla  \vh) - (\kappa^2 \ExtLoch{j} \tauh, \vh)_{\Oj} &= 0  , \quad\text{for all } \vh \in  \VFEMO(\Oj)   , \\
    (\ExtLoch{j} \tauh) |_{\partial \Oj} &= \tauh   .
\end{aligned}
\end{gather}
Let us introduce the discrete extension $\ExtGh: \VFEM(\Gamma) \to \VFEM(\Omega)$, which first restricts functions to the boundaries of the subdomains $\Oj$ and then applies $\ExtLoch{j}$. 

The discrete ACMS space, indicated with the subscript $A$, is then defined on the domain decomposition $\D$ as the span of the extended discrete edge modes and discrete vertex basis functions:
\begin{align*}
    \VACMSh  &= \mbox{ span} \{ \ExtGh \tauhe, \ \forall e \in \E, \ i=1,...\Ie \} \ \cup \mbox{ span} \{ \ExtGh \phiPh, \ \forall q \in \V \}  .
\end{align*}
The weak formulation~\eqref{eq:ACMS_Weak_Formulation} can be rewritten in the ACMS space: find $\uACMSh \in \VACMSh $ such that
\begin{align} \label{eq:ACMS_Weak_Formulation_Discrete}
\sesH(\uACMSh, \vACMSh) = G(\vACMSh), \quad \text{for all } \vACMSh \in \VACMSh   .
\end{align}
Finally, the weak formulation in the ACMS space~\eqref{eq:ACMS_Weak_Formulation_Discrete} turns into
\begin{align} \label{eq:ACMS_System}
    \SACMS \uvecACMS = \gvecACMS  ,
\end{align}
where $\uvecACMS \in \RR^{\NACMS}$ is the solution vector, $\SACMS \in \RR^{\NACMS \times \NACMS}$ is the system matrix and $\gvecACMS \in \RR^\NACMS$ is the right-hand side vector. Here, $\NACMS = \dim \VACMSh$ is the space dimension.

\subsection{Computational Costs}  \label{sec:computational_cost}
In the following, we give some insights on the expected computational costs, in terms of FLOPS, when using the discrete ACMS method. For more details, we refer to \cite{Heinlein2015}, where an efficient and parallel implementation for $p=1$ is discussed. We are interested in the dependence of the cost on the number of modes $\Ie$, the number of subdomains $J$ and related quantities, such as $|\V|, |\E|$, and the corresponding dimension of the ACMS space $N_A = |\V| + |\E| \Ie$. 

To ease the presentation, we make the following assumptions:
\begin{itemize}
    \item all subdomains $\Oj \in \D$ of the domain decomposition are of comparable size, and correspondingly also the dimensions $\NFEMj = \dim \VFEM(\Oj)$ are similar,
    \item all edges $e \in \E $ are of comparable size, and correspondingly also the dimensions $\NFEMe = \dim \VFEM(e)$ are similar.
\end{itemize}
For a given mesh size $h$ and polynomial order $p$ of the underlying finite element space, we can introduce the following two costs
\begin{align*}
    &\Cext:\text{ cost of solving the extension problem \eqref{eq:localHarmonicExtensionDiscrete} on a subdomain $\Oj \in \D$},\\
    &\Ceig:\text{ cost of solving the eigenvalue problem \eqref{eq:edge_mode_discrete} on an edge $e \in \E$}.
\end{align*}

\paragraph{Edge basis computation}
The cost for the computation of the edge modes is given by solving the eigenvalue problem \eqref{eq:edge_mode_discrete} for all edges $e \in \E$ and by the application of the extension operator to the corresponding neighboring subdomains for each mode (and each edge) separately; see \Cref{rem:extension}. Thus, it is given by the respective contributions
\begin{align} \label{eq:edge_cost}
    \O \left(|\E| \Ceig \right) + \O( |\E| \Ie \Cext).
\end{align}

\begin{remark} \label{rem:eigenvalue}
    When the eigenvalue computation has to be done for each edge separately, the first term in \eqref{eq:edge_cost} might dominate when $\Ie$ grows. Although, in practice, we expect that the problem structure can be exploited such that the eigenvalue problem only needs to be solved $\O(1)$ times; see \Cref{Sec:Timings}.
\end{remark}

\paragraph{Vertex basis computation}
The costs of the vertex basis functions, where we only consider the application of the harmonic extensions, is given by
\begin{align} \label{eq:vertex_cost}
    \O\left( |\V| \Cext\right).
\end{align}
Note that, for the above cost, we assumed that the restrictions on the edges of functions $\phiPh$ which fulfill property \eqref{eq:vertex_harmonic} are either given or very cheap to compute. This assumption is motivated by the fact that, in case the edges of the domain decomposition are straight lines, the restriction of $\phiPh$ on an edge $e \in \E$ is given by a linear function.

\paragraph{Assembly}
From an implementation perspective, the ACMS method exhibits some similarities to a high-order finite element method, which can be utilized during setup and assembly; see also \Cref{rem:hp-acms,rem:parallel}. 
More precisely, as already discussed in \cite{Hetmaniuk2010,Heinlein2015}, the assembly is done in the usual way, via a loop over the subdomains $\Oj \in \D$, see also, e.g., \cite{Braess2007}. 
On each cell $\Oj$, one then first computes the local ACMS basis functions (as discussed in the previous paragraph) and then calculates the local ACMS stiffness matrix $\SACMS$ by means of the stiffness matrix $\SFEM$ of the underlying finite element space $\VFEM(\Oj)$. With a slight abuse of notation let us write the local contributions on $\Oj \in \D$ as 
\begin{align}\label{eq::SACMS_via_trasf_loc}
    \SACMS|_{\Oj} = \BB_j^T \SFEM|_{\Oj} \BB_j, 
\end{align}
where the restriction is understood as taking only into account degrees of freedom of basis functions (of the corresponding finite element spaces) which have a non-zero support on $\Oj$. Likewise, we indicate with $\BB_j$ the matrix where each column corresponds to a restriction of the ACMS-basis functions supported on $\Oj$. Consequently, the column dimension of $\BB_j$ scales with $\O(\Ie)$. Since $\SFEM|_{\Oj}$ is sparse with $\O(\NFEMj)$ bandwidth, the costs for \eqref{eq::SACMS_via_trasf_loc} are $\O(\NFEMj \Ie^2)$. Since this has to be done on each subdomain the total cost for the assembly are given by 
\begin{align*}
    \O\left(J \NFEMj \Ie^2 \right).
\end{align*}

\paragraph{Solving the system}
Finally, the cost for solving the system \eqref{eq:ACMS_System} is $\O(\NACMS^3)$ in case the matrix $\SACMS$ is dense.
However, since the matrix $\SACMS$ shows a sparsity pattern, due to the aforementioned local support of the ACMS basis functions, as depicted in \cref{figsparsity}, these costs can be further reduced to 
\begin{align} \label{eq:costs_solve_bandwidth}
    \O(N_A b_A^2), 
\end{align}
where $b_A$ is the bandwidth of $\SACMS$; see \cite[Ch. 2.1.1]{Quarteroni2008NumericalEquations}. We note that on structured grids this estimate might be pessimistic; see \cite{LiptonRoseTarjan1979}. Indeed, we observe a better scaling in the following section.

\begin{remark} \label{rem:extension}
    The cost for the extension $\Cext$ on a given domain $\Oj \in \D$,  corresponds to solving a sparse system associated to problem \eqref{eq:localHarmonicExtensionDiscrete} and is of order $\O(\NFEMj \NFEMe^2)$; see~\cite[Ch.~2.1.1]{Quarteroni2008NumericalEquations}. 
    In the implementation of the extension operators, we only compute the factorization once (using a sparse Cholesky solver), and subsequently apply the forward/backward substitution for each basis function independently. 
    Therefore, the contribution of $\Cext$ is split into the factorization and application costs. The factorization is reused for both edge and vertex basis functions and, therefore, is carried out once per subdomain, accounting for a total cost of $\O( J \NFEMe^4)$. The application is carried out $\Ie$ times per subdomain in the edge basis computation \eqref{eq:edge_cost}, meaning $\O( |\E| \Ie \NFEMe^3)$, and only once per subdomain in the vertex basis computation \eqref{eq:vertex_cost}, which is $\O( |\E| \NFEMe^3)$.  
    
\end{remark}

\begin{remark}\label{rem:eigenvalue_computation}
    We emphasize that $\Ceig$ clearly depends on the choice of the number of edge modes $\Ie$. In fact, the cost for solving the eigenvalue problem on an edge $e \in \E$ is of order $\O(\NFEMe^2 \Ie)$, using, e.g., an Arnoldi GMRES, see \cite[Ch. 2.6.2]{Quarteroni2008NumericalEquations}. 
\end{remark}

\begin{remark}\label{rem:eigenvalue_dimension}
    We remind that, since it always holds that $\Ie \le \NFEMe$, the computation of a high number of edge modes will require enriching the underlying finite element space $\VFEMO(e)$. That can be obtained by either refining the mesh size or increasing the polynomial degree of approximation.
\end{remark}

\begin{remark} \label{rem:hp-acms}
If the ACMS method is interpreted as a high-order finite element method, increasing the maximum number of modes $\Ie$ is similar to having a $p$-FEM method, while increasing the number of domains $J$ can be related to an $h$-FEM method. 
\end{remark}

\begin{remark} \label{rem:parallel}
    So far, we have not considered the advantages of potential parallel implementation. Nevertheless, we want to emphasize that, as with standard finite element methods, there is considerable room for improvement. For example, the assembly process and the independent application of the extensions (see \Cref{rem:extension}) can be parallelized.
\end{remark}

\subsection{Comparison of costs and timings for a prototype crystal} \label{Sec:Timings}

\begin{figure}[ht]
    \centering
    \begin{tikzpicture}
        \node[] at (0,0) {\includegraphics[width=6cm, clip = true, trim = 2.5cm 0 2.5cm 0]{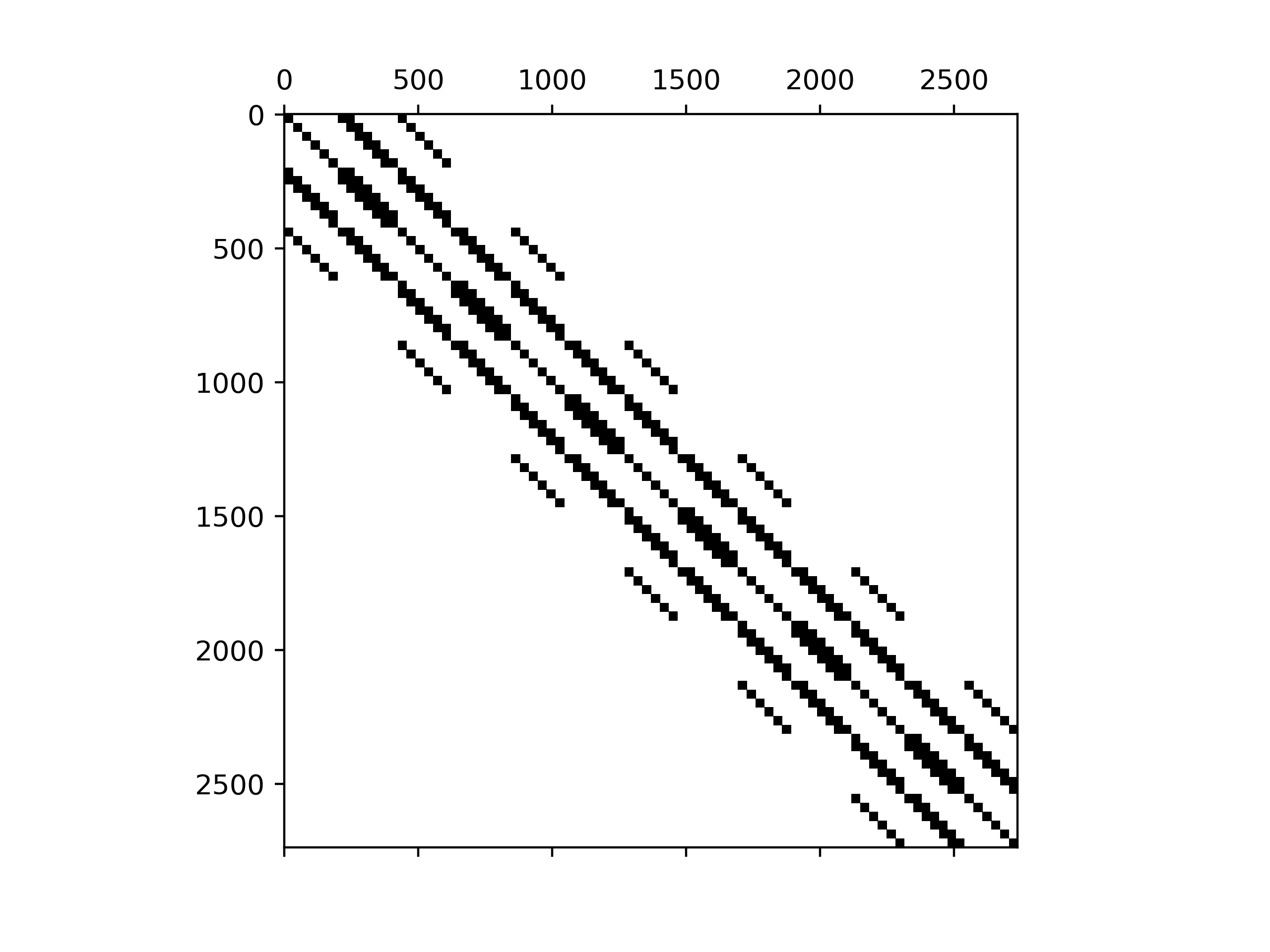}};
        \begin{scope}[shift = {(4,-2)}]
            \input{dofnumbering.tex}
        \end{scope}
    \end{tikzpicture}
    \caption{Left: Sparsity pattern of $\SACMS$ for $\Ie = 32$ and $J = 36$ domains. The bandwidth is approximately $b_A \sim 3\sqrt{J} \Ie$, hence $b_A \sim 384$. Right: dof numbering of a crystal with $J = 4$ and $\Ie = 2$.}
    \label{figsparsity}
\end{figure}

In this section, we want to discuss the above derived computational costs and compare them to the timings of our implementation.
We choose $a = \beta = \kappa = 1$. 
As the computational domain we consider a square-shaped prototype photonic crystal structure $\Omega = [0,\sqrt{J}] \times [0, \sqrt{J}]$ decomposed into $J$ unit square subdomains that we refer to as unit cells. 
For a given mesh size $h$ and order $p$, considering a quasi uniform and shape regular triangulation of $\Omega$, we have $\NFEMj = \NFEMe^2 \sim (p/h)^2$ for every subdomain. It is clear that, for this setting, the assumptions from \Cref{sec:computational_cost} are fulfilled. Additionally, we consider a mesh such that vertical and horizontal edges have the same (one dimensional) mesh and thus the eigenvalue problem \eqref{eq:localHarmonicExtensionDiscrete} has to be solved only once (and the cost is ignored); see \Cref{rem:eigenvalue}. Since $|\E|\sim J$ and $|\V| \sim J$, and hence $\NACMS = |\E| \Ie + |\V| \sim J \Ie$, the total computational costs are given by 
\begin{align} \label{eq:acms_total_costs}
    \underbrace{\O(J \Ie \Cext)
    }_{\text{basis computation}}
    +
    \underbrace{\O(J \NFEMj \Ie^2 )
    }_{\text{assembly}} 
    +
    \underbrace{\O( J^2 \Ie^3)
    }_{\text{solve}},
\end{align}
where we used that the bandwidth of our structured prototype crystal is $b_A \sim 3 \sqrt{J} \Ie$. The latter follows by considering the following degrees of freedom (dof) numbering of the ACMS basis: we start with the lower left vertex and number all dofs on vertices and edges increasingly along the $x$-axis, i.e., we consider only horizontal edges. Next, all dofs on vertical edges are numbered increasingly along the $y$-axis of the lowest row of the domain. This pattern is repeated for each row.
Then the aforementioned bandwidth follows since we have $\sqrt{J}$ unit cells in each direction and $\Ie$ basis functions per edge, see also \cref{figsparsity} for an illustration of the sparsity pattern and the dof numbering.
For the total costs \eqref{eq:acms_total_costs}, we observe that, similarly to standard finite element methods, solving the system could potentially be the dominant factor for very large $J$ and $\Ie$. However, for the discrete ACMS method, we expect that $N_A$ is small enough that the assembly and the basis computation are the limiting factors in the computations. This is motivated by the limit case when $\Ie = \NFEMe = \sqrt{\NFEMj}$, then, taking into account \Cref{rem:extension,rem:eigenvalue_computation}, the full cost \eqref{eq:acms_total_costs} becomes 
\begin{align*} \label{eq:acms_total_costs_limitcase}
    \underbrace{\O(J \Ie^4)
    }_{\text{basis computation}}
    +
    \underbrace{\O(J \Ie^4 )
    }_{\text{assembly}} 
    +
    \underbrace{\O( J^2 \Ie^3)
    }_{\text{solve}}.
\end{align*}
When doing $hp$--refinements while keeping the domain decomposition fixed, the number of modes $\Ie$ becomes much larger than the number of subdomains $J$ and the cost is dominated by the basis computation and assembly.

\tikzstyle{vertex}=[scale=0.7,fill=blue!20] 
\tikzstyle{vertexedge}=[scale=0.7,fill=red!20] 
\tikzstyle{edge}=[black]

\begin{figure}[ht]
\centering
\includegraphics[scale=1]{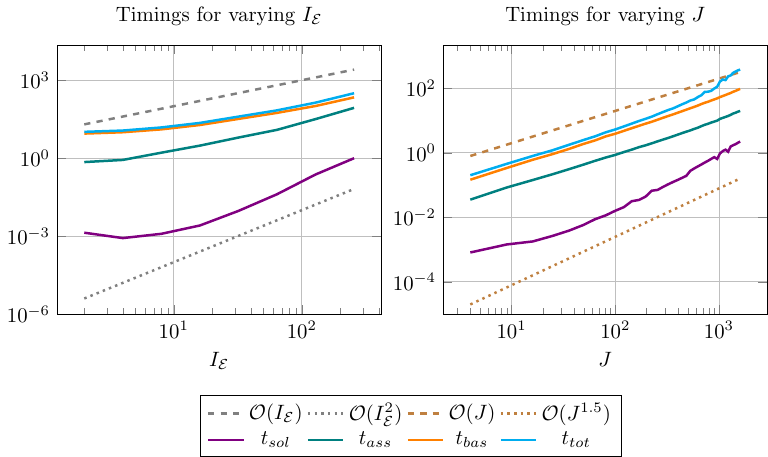}
\caption{Timings and reference lines for varying number of edge modes $\Ie$ (left) and number of sub domains $J$ (right).}
    \label{fig:scaling}
\end{figure}

\paragraph{Scaling with respect to $\Ie$} 
In \Cref{fig:scaling} (left), we present the timings for the computation of the basis functions $t_{bas}$, the assembly $t_{ass}$, the time to solve the system $t_{sol}$, and the sum of these timings $t_{tot}$. For our computation, we choose $h = 0.01$, $p = 6$, $J = 4$, and we set $\Ie = 2^l$ with $l = 1,\ldots,9$. According to \eqref{eq:acms_total_costs}, we expect a linear scaling $\O(\Ie)$ for the basis computation, a quadratic scaling $\O(\Ie^2)$ for the assembly and a cubic scaling $\O(\Ie^3)$ for solving the system. While the times of the basis computations and the assembly given in \Cref{fig:scaling} (left) are in accordance with the theory, we see a  better scaling for the solving times. Note, that \eqref{eq:costs_solve_bandwidth} might be pessimistic; see \cite{LiptonRoseTarjan1979}. In our computations we used SciPy to solve the system, which uses the sparse UMFPACK solver; see \cite{UMFPACK2004,SciPy2020}.
In practice, considering such a high number of modes might not be necessary due to the fast error reduction; see \eqref{eq:L2_error_bound_kappa_dep}.

\paragraph{Scaling with respect to $J$} 
In \Cref{fig:scaling} (right), we again present the timings for the computation of the basis functions $t_{bas}$, the assembly $t_{ass}$, the time to solve the system $t_{sol}$, and the sum of these timings $t_{tot}$. For our computation, we choose $h = 0.1$, $p = 4$, $\Ie = 4$, and we set $J = l^2$ with $l = 1,\ldots,40$. According to \eqref{eq:acms_total_costs} we expect a linear scaling $\O(J)$ for the basis computation and the assembly, and a quadratic scaling $\O(J^2)$ for solving the system. Again, while the timings for the basis computations and the assembly are in accordance to the theory, the solver provided by SciPy now scales with the optimal rate $\O(J^{3/2})$; see  \cite{LiptonRoseTarjan1979}.

\section{Numerical results} \label{Sec:NumericalResults}

The main goal of this section is to exemplify how to choose the different parameters of the discrete ACMS method in order to solve large scale problems. 
Namely, we will present a study on suitable choices of the mesh size of the (underlying) finite element space, the polynomial order, and the number of edge modes needed to achieve a certain accuracy, depending on the material parameters and the chosen wavenumber. The computations were done using the finite element library Netgen/NGSolve, see \cite{NGSolve}. 
Further, for reproduction purposes, all data can be found in the repository \cite{NGSolveACMS}.

We consider three different numerical examples. First, in \Cref{test:numex1}, we investigate the effect of using an underlying higher order FEM for the ACMS basis computation on the method accuracy in relation to the wavenumber. 
Then, in \Cref{test:numex2}, we test a configuration with discontinuous coefficients and, given the flexibility of the ACMS method, we examine how the chosen domain decomposition impacts the accuracy and efficiency of the method. 
Finally, in~\Cref{test:numex3}, we show the applicability and flexibility of the method on a large but finite-size structure that has defects.

\subsection{Example 1: Circular domain} \label{test:numex1}
\begin{figure}[ht]
\centering
\begin{subfigure}[b]{0.45\textwidth}
\centering
\begin{tikzpicture}[>=latex, scale = 2.4]
    \def\r{0.03}
    \def\R{1}
    \def\l{0.1mm}
    \fill[white, fill opacity=0.05, draw = black, line width = \l] (0,0) circle (\R);

    \begin{scope}[densely dotted,ultra thick]
        
    \draw[black, line width = \l]  (\R, 0)  to  (0,\R) ;
    \draw[black, line width = \l] (-\R, 0) to   (0,\R);
    \draw[black, line width = \l] (-\R, 0) to  (0,-\R);
    \draw[black, line width = \l]  (0,-\R) to    (0,\R);

    \draw[black, line width = \l] (-\R,0) to (\R,0);
    \draw[black, line width = \l] (0,-\R) to (\R,0);
 
    \fill[black, line width = \l] (0,0) circle     (\r);
    \fill[black, line width = \l] (0,\R) circle    (\r);
    \fill[black, line width = \l] (\R,0) circle    (\r);
    \fill[black, line width = \l] (-\R,0) circle   (\r);
    \fill[black, line width = \l] (0,-\R) circle  (\r);
    \end{scope}
\end{tikzpicture}
\caption{ \label{fig:CircularDomain}}
\end{subfigure}
\hfill
\begin{subfigure}[b]{0.45\textwidth}
\centering
\begin{tikzpicture}[>=latex, scale = 2.4]
    \def\r{0.03}
    \def\R{1}
    \def\l{0.3mm}
    \fill[white, draw = black, line width = \l] (0,0) circle (\R);
    \begin{scope}[draw = gray]
    \input{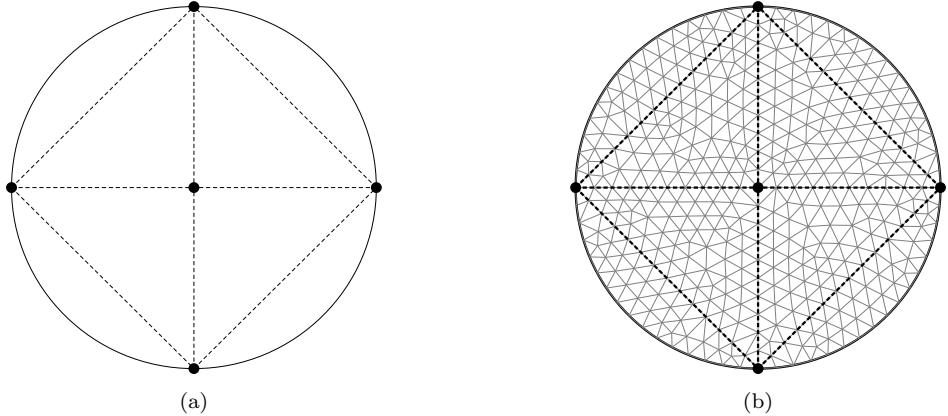}
    \end{scope}
    \begin{scope}[densely dotted,ultra thick]
    \draw[black, line width = \l]  (\R, 0)  to  (0,\R) ;
    \draw[black, line width = \l] (-\R, 0) to   (0,\R);
    \draw[black, line width = \l] (-\R, 0) to  (0,-\R);
    \draw[black, line width = \l]  (0,-\R) to    (0,\R);

    \draw[black, line width = \l] (-\R,0) to (\R,0);
    \draw[black, line width = \l] (0,-\R) to (\R,0);
 
    \fill[black, line width = \l] (0,0) circle     (\r);
    \fill[black, line width = \l] (0,\R) circle    (\r);
    \fill[black, line width = \l] (\R,0) circle    (\r);
    \fill[black, line width = \l] (-\R,0) circle   (\r);
    \fill[black, line width = \l] (0,-\R) circle  (\r);
    \end{scope}
\end{tikzpicture}
\caption{\label{fig:CircularDomain_Mesh}}
\end{subfigure}    
\caption{In \subref{fig:CircularDomain}, the domain $\Omega$ and the domain decomposition $\D$. In \subref{fig:CircularDomain_Mesh}, the corresponding coarse finite element mesh with $h=0.1$. The vertices in $\V$ are marked with dots.}
\end{figure}

We consider~\eqref{eq:model_strong} on the unit disc $\Omega=B_1(0)$ and we set the constant functions $a  = 1$,  $\beta = 1$, $c =1$, and $f =0$. We define the boundary source term $g$ such that the problem admits the plane wave $u  = e^{- \cplu \kappa (0.6x+0.8y)}$ as its analytical solution, with constant wavenumber $\kappa=\omega$. We mention that this setting is as in \cite[Section~5.1]{GiammatteoHeinleinSchlottbom2024}.
We choose a domain decomposition $\D$ as shown in \Cref{fig:CircularDomain}, with $J = 8$ subdomains, $12$ edges in $\E$, and $5$ vertices in $\V$. 
In ~\Cref{fig:CircularDomain_Mesh}, we show a coarse triangulation that is uniformly refined in order to compute the approximated ACMS solution. Below, we will compare the relative $L^2$--error, i.e.,  
$\| u - \uACMSh \|_{\LO}/\| u\|_{\LO}$
of the approximated ACMS solution against the analytical solution; this will verify the correct implementation of the high-order ACMS method. 
In this example, we examine the effects of increasing the maximum number of computable edge modes by either refining the mesh size or driving up the polynomial order of approximation; in both cases, $\NFEMe$ is increased, as commented in~\Cref{rem:eigenvalue_dimension}.

\begin{figure}
\centering \small 
\setlength\tabcolsep{0.55em}
\begin{tikzpicture} 

    \begin{groupplot}[group style ={group size = 3 by 1, horizontal sep=0.5cm},
        width = 0.39\textwidth, height = 0.4\textwidth,
    ymode=log,
    legend style={
      legend columns=4,
      at={(0.5,-0.4)},
      anchor=south,
      },
    xlabel = \( p \), 
    ymin = 5.0e-7, 
    xmajorgrids=true, ymajorgrids=true,
    xtick={1,2,3,4,5},
    xticklabels={1,2,3,4,5},
    ytick={1e0, 1e-1, 1e-2, 1e-3, 1e-4, 1e-5, 1e-6},
    yticklabels={$1$, $10^{-1}$, $10^{-2}$, $10^{-3}$, $10^{-4}$, $10^{-5}$, $10^{-6}$},
    scaled y ticks = false,
    cycle list name=philipcolors,
    ]
   \pgfplotsset{cycle list shift=8}
   
           \nextgroupplot[ xlabel =  {$h = 0.025, \Ie = 32$}, ]
           
        \addplot
        coordinates { (1, 8.204e-2) (2, 3.957e-4) (3, 2.382e-4) (4, 2.379e-4) (5, 2.379e-4) }; 
        
        \addplot
        coordinates {(1, 6.041e-1) (2, 5.967e-3) (3, 1.430e-3) (4, 1.422e-3) (5, 1.422e-3)    }; 
            
        \addplot
        coordinates {(1, 1.278e0) (2, 1.622e-1) (3, 3.991e-2) (4, 3.963e-2) (5, 3.963e-2)  }; 

        \addplot
        coordinates {(1, 1.217e0) (2, 1.296e0) (3, 9.151e-1) (4, 8.559e-1) (5, 8.554e-1)};  

   \nextgroupplot[ xlabel = {$h = 0.0125, \Ie = 64$}, yticklabels=\empty ]
   
    \addplot
    coordinates { (1, 2.019e-2) (2, 4.344e-5) (3, 2.940e-5) (4, 2.937e-5) (5, 2.937e-5) };  
    
    \addplot
    coordinates {(1, 1.598e-1) (2, 4.199e-4) (3, 1.330e-4) (4, 1.328e-4) (5, 1.328e-4) };  
        
    \addplot
    coordinates {(1, 1.074e0) (2, 1.010e-2) (3, 7.941e-4) (4, 7.847e-4) (5, 7.846e-4) }; 

    \addplot
    coordinates {(1, 1.335e0) (2, 2.890e-1) (3, 3.195e-2) (4, 3.138e-2) (5, 3.138e-2) };  

\nextgroupplot[xlabel = {$h = 0.00625, \Ie = 128$}, yticklabels=\empty,
    legend style={
      legend columns=4,
      at={(-0.6,-0.3)},
      anchor=north
  }
    ]
    
    \addplot
        coordinates {(1, 5.116e-3) (2, 5.359e-6) (3, 3.699e-6) (4, 3.695e-6) (5, 3.695e-6)};  
        
        \addplot
        coordinates {(1, 4.097e-2) (2, 3.899e-5) (3, 1.540e-5) (4, 1.538e-5) (5, 1.538e-5) };  
            
        \addplot
        coordinates {(1, 3.209e-1) (2, 6.860e-4) (3, 6.897e-5) (4, 6.867e-5) (5, 6.867e-5) };  

        \addplot
        coordinates {(1, 1.497e0) (2, 2.009e-2) (3, 4.899e-4) (4, 4.731e-4) (5, 4.730e-4)};  

        \legend{$\kappa = 16$, $\kappa = 32$, $\kappa = 64$,$\kappa = 128$ }
    \end{groupplot}
    
    \node[ellipse,
    draw = blue,
    fill = green!20!blue,
    opacity=0.2,
    minimum width = 0.4cm, 
    minimum height = 1.6cm] (e) at (11.74, 0.97) {};
    
\end{tikzpicture}
   
\caption{Example in \Cref{test:numex1}, on a circular domain with decomposition $\D$ depicted in \Cref{fig:CircularDomain}. The $L^2$--relative error $\| u - \uACMSh \|_{\LO}/\| u\|_{\LO}$ is computed against the exact solution for a fixed mesh size and number of edge modes. We choose $h = 0.025$ and $\Ie = 32$ (left), $h = 0.0125$ and $\Ie = 64$ (middle), $h=0.00625$ and $\Ie = 128$ (right). Test with increasing degree of approximation $p=1,2,3,4,5$ (horizontal axis), for different wavenumber values $\kappa = 16 , \ 32, \ 64 , \ 128 $. In the blue region, the values are highlighted for a better comparison with \Cref{fig:Circle_Convergence_Full}.
\label{fig:Circle_Convergence_Saturation}}
\end{figure}
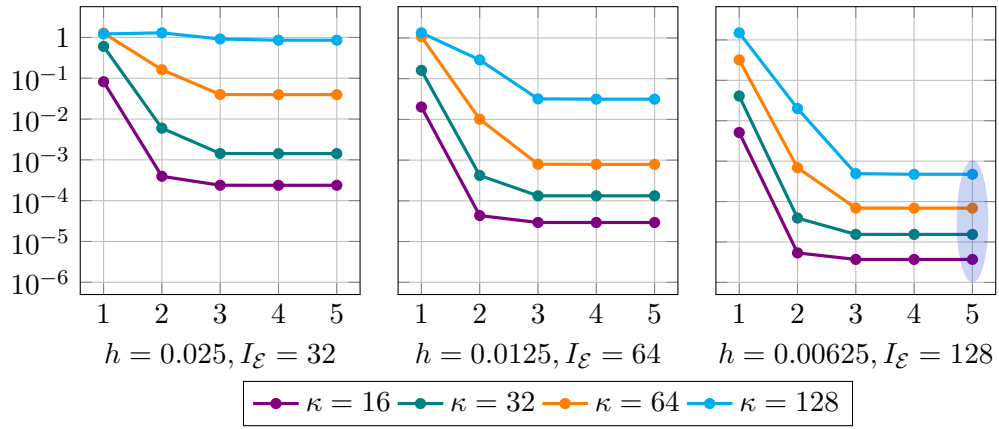

First, we aim to investigate how the degree of approximation of the underlying finite element space, used to approximate the computed ACMS basis functions, affects the accuracy of the ACMS method.
In \Cref{fig:Circle_Convergence_Saturation}, we show the $L^2$--relative errors for different discretization parameters: $h = 0.025$ and $\Ie = 32$ (left), $h = 0.0125$ and $\Ie = 64$ (middle), and $h=0.00625$ and $\Ie = 128$ (right). We consider an increasing degree of approximation $p=1,2,3,4,5$, and we repeat the tests for different wavenumbers $\kappa = 16, 32, 64, 128$.
On the one hand, we observe that the approximation gets worse for higher wavenumbers. In particular, there is hardly any error reduction for $p=1$ and $\kappa=128 \ (\circ)$, which might be explained by a pollution effect. 
On the other hand, the error decreases approximately by an order of magnitude if $\Ie$ is doubled for each $p\geq 3$ fixed.
The observed saturation effect in \Cref{fig:Circle_Convergence_Saturation} occurs, because of the fixed number of employed edge modes and the decreasing value for $h$.

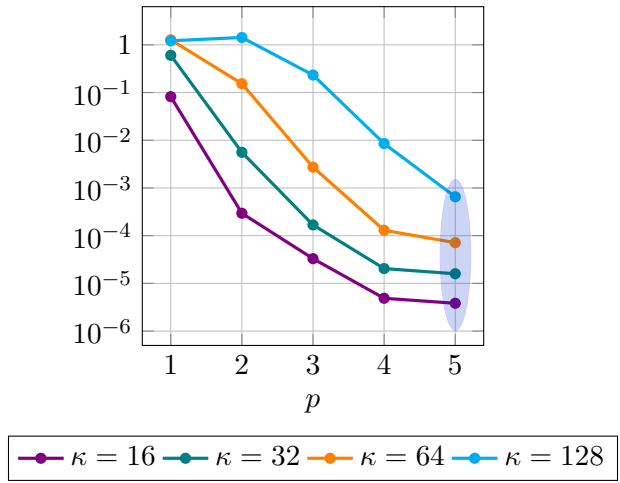
\begin{figure}
\centering\small \setlength\tabcolsep{0.55em}
\begin{tikzpicture}
    \begin{axis}[    
    width=0.45\textwidth,height=0.45\textwidth,
    ymode=log,
    legend style={
      legend columns=4,
      at={(0.5,-0.4)},
      anchor=south,
      },
    xlabel = \( p \), 
    ymin = 5.0e-7, 
    xmajorgrids=true, ymajorgrids=true,
    xtick={1,2,3,4,5},
    xticklabels={1,2,3,4,5},
    ytick={1e0, 1e-1, 1e-2, 1e-3, 1e-4, 1e-5, 1e-6},
    yticklabels={$1$, $10^{-1}$, $10^{-2}$, $10^{-3}$, $10^{-4}$, $10^{-5}$, $10^{-6}$},
    scaled y ticks = false,
    cycle list name=philipcolors,
    ]
   \pgfplotsset{cycle list shift=8}

    \addplot
        coordinates {(1, 8.204e-2) (2, 2.949e-4) (3, 3.296e-5) (4, 4.869e-6) (5, 3.827e-6)};  
        
        \addplot
        coordinates {(1, 6.041e-1) (2, 5.583e-3) (3, 1.681e-4) (4, 2.043e-5) (5, 1.589e-5)};  
            
        \addplot
        coordinates {(1, 1.278e0) (2, 1.532e-1) (3, 2.735e-3) (4, 1.297e-4) (5, 7.127e-5)};  

        \addplot
        coordinates {(1, 1.217e0) (2, 1.428e0) (3, 2.331e-1) (4, 8.528e-3) (5, 6.533e-4)};  
        
    \legend{$\kappa = 16$, $\kappa = 32$, $\kappa = 64$,$\kappa = 128$ }
    \end{axis}
    \node[ellipse,
    draw = blue,
    fill = green!20!blue,
    opacity=0.2,
    minimum width = 0.4cm, 
    minimum height = 2cm] (e) at (4.12,1.2) {};
\end{tikzpicture}
\caption{Example in \Cref{test:numex1}, on a circular domain with decomposition $\D$ depicted in \Cref{fig:CircularDomain}. The $L^2$--relative error $\| u - \uACMSh \|_{\LO}/\| u\|_{\LO}$ is computed against the exact solution for mesh size $h = 0.025$ and for different wavenumber values $\kappa = 16 , \ 32 , \ 64 , \ 128 $. The modes considered are $\Ie = 32$ for degree $p=1$, $\Ie = 64$ for degree $p=2,3$ and $\Ie = 128$ for degree $p=4,5$.
In the blue region, the values are highlighted for a better comparison with \Cref{fig:Circle_Convergence_Saturation}.
\label{fig:Circle_Convergence_Full}}
\end{figure}

To overcome this saturation effect, it is necessary to employ more edge modes while keeping $h$ fixed. 
In~\Cref{fig:Circle_Convergence_Full}, we adjust the number of modes based on the degree of approximation, since, as noted in \Cref{rem:eigenvalue_dimension}, the maximum number of possible edge modes $\Ie$ increases for a higher polynomial order. 
We choose mesh size $h=0.025$ and the modes considered are $\Ie = 32$ for degree $p=1$, $\Ie = 64$ for degree $p=2,3$ and $\Ie = 128$ for degree $p=4,5$.
This choice was made because these numbers of modes are the largest values we can compute for each degree of approximation while still being powers of two, ensuring consistency with the previous tests.
As expected, we can now see convergence without saturation and correspondingly the full benefit of choosing a high-order approximation. The edge modes are better approximated because the smooth eigenfunctions, which underly the ACMS basis functions, are effectively approximated using a $p$-method. Additionally, the number of computable edge modes increases, which enriches the approximation space of the ACMS method and further reduces the approximation error.
To further highlight this observation, we point out that we obtain the same error accuracy if we compare the error in~\Cref{fig:Circle_Convergence_Saturation} (right) (where $h = 0.00625$)  and in~\Cref{fig:Circle_Convergence_Full} (where $h = 0.025$), for example, for degree $p=5$ and $\Ie = 128$ (errors within the blue ellipses). This means that, already for the coarsest mesh, we obtain the most accurate solution if we use high polynomial order $p$ and employ a suitable number of modes. 

Before proceeding to the next example, we want to verify numerically the dependence between the number of edge modes and the wavenumber that we observed in our theoretical result \eqref{eq:L2_error_bound_kappa_dep}. 
In~\Cref{fig:Circle_Wavenumber}, we plot in logarithmic scale the relative $L^2$--error for varying number of edge modes $\Ie = 1, ..., 128$ for different wavenumber $\kappa = 16, 32, 64, 128$ and for fixed mesh size $h=0.025$ and degree of approximation $p=5$.
Qualitatively, the curves show the same behavior.
First the error stagnates, and then there is a significant drop after which the asymptotic convergence rate is achieved.
Stagnation of error might be explained by requiring a minimal resolution condition, compare this to the situation of piecewise linear approximation, where $h\kappa\approx 1/10$ is required to be able to interpolate a corresponding wave accurately on the given grid.
The transition from stagnating errors to monotonically decreasing errors is marked with a black circle.
We observe that the required number of modes to enter this monotone behavior almost doubles if the wavenumber is doubled. 
This might be related to a minimal resolution condition.
Moreover, for this example, we may state that the adjoint approximability constant $\sigma^*$ behaves better than the scaling $\O(\Cstab\|\kappa\|_\infty^2 / \Ie)$ predicted by theory, see the discussion following \eqref{eq:ACMS_Weak_Formulation}.
We conclude that, in this example, when doubling the wavenumber, the errors achieve the predicted convergence rates and enter the asymptotic regime when the number of edge modes is also doubled.

\pgfplotstableread[]{circle_omega_1.0_16.0_32.0_64.0_128.0_maxH_0.025_order_5_EM_1_128__table.out}\WavenumberSweepCircle

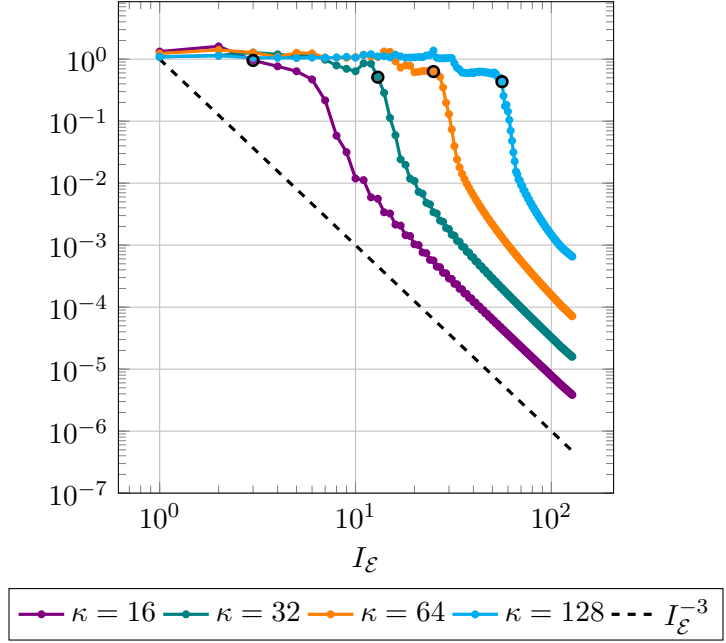
\begin{figure}[ht]
\centering\small \setlength\tabcolsep{0.55em}
\begin{tikzpicture}
  \begin{axis} [ 
    width=0.6\textwidth,height=0.6\textwidth,
    ymode=log,
    xmode=log,
    legend style={
      legend columns=5,
      at={(0.5,-0.3)},
      anchor=south,
      },
    xlabel = \( \Ie \), 
    ymin = 1.0e-7, 
    xmajorgrids=true, ymajorgrids=true,
    scaled y ticks = false,
    cycle list name=philipcolors,
    ]

\pgfplotsset{cycle list shift=15}

   \addplot table[x=Ie,y expr={\thisrowno{2}} ]{\WavenumberSweepCircle}; 
    \addplot table[x=Ie,y expr={\thisrowno{3}} ]{\WavenumberSweepCircle}; 
    \addplot table[x=Ie,y expr={\thisrowno{4}} ]{\WavenumberSweepCircle}; 
    \addplot table[x=Ie,y expr={\thisrowno{5}} ]{\WavenumberSweepCircle}; 

    \addplot[black,dashed,very thick,mark size=3pt] 
    table[x=Ie,y expr={(\thisrowno{0})^(-3)} ]{\WavenumberSweepCircle};

   \addplot[  color=black, mark=*, fill = white, fill opacity=0.1, mark size = 2pt , mark options={line width=1pt} ]  coordinates {(3, 9.53968920e-01)  };
   
    \addplot[ color=black,  mark=*, fill = white, fill opacity=0.1,  mark size = 2pt , mark options={line width=1pt}] coordinates {(13, 5.12090439e-01) };
    \addplot[ color=black,  mark=*, fill = white, fill opacity=0.1,  mark size = 2pt , mark options={line width=1pt}] coordinates {(25, 0.6309664202496508 )};
    \addplot[ color=black,mark=*, fill = white, fill opacity=0.1,  mark size = 2pt, mark options={line width=1pt}]   coordinates {(56,  0.43286103078918325)   };
 
      \legend{$\kappa = 16$, $\kappa = 32$, $\kappa = 64$,$\kappa = 128$, $\Ie^{-3}$}
    \end{axis}  
\end{tikzpicture}
\caption{Example in \Cref{test:numex1}, on a circular domain with decomposition $\D$ depicted in \Cref{fig:CircularDomain}. The $L^2$--relative error $\| u - \uACMSh \|_{\LO}/\| u\|_{\LO}$, shown in logarithmic scale, is computed against the exact solution for a fixed order of approximation $p = 5$ and with mesh size $h= 0.025$. We consider an increasing number of edge modes per edge $\Ie = 1, ...,128$ for different wavenumber values $\kappa = 16 , \ 32, \ 64 , \ 128 $. 
The small black circles indicate the index $\Ie$ after which the error decays monotonically.
\label{fig:Circle_Wavenumber}}
\end{figure}

\subsection{Example 2: Square crystal with circular pores} \label{test:numex2}
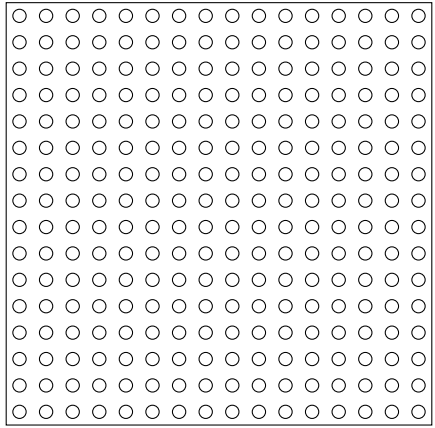
\begin{figure}[ht]
\centering
\begin{tikzpicture}[>=latex, scale = 0.35]
  \def\Ncell{16}
  \def\incl{2}
  \pgfmathsetmacro\I{0.5/\incl}
  \pgfmathsetmacro\N{\Ncell/\incl}
  \def\r{0.25}

  \fill[white, fill opacity=0.2, draw = black] (0,0) rectangle (\Ncell,\Ncell);
  
  \foreach \i in {1, ..., \Ncell} {
    \foreach \j in {1, ..., \Ncell} {
          \draw[] (-0.5+\i,-0.5+\j) circle (\r);
    }
  }
\end{tikzpicture}
\caption{$\Omega=(0,16)\times(0,16)$ with   $16 \times 16$ square unit cells. In this case the decomposition $\Dincl{16}$ contains only one subdomain $\Omega$.\label{fig:crystal_domain}} 
\end{figure}

We consider~\eqref{eq:model_strong} on a square domain $\Omega=[0,16] \times [0,16]$ that is meant to model a finite-size photonic crystal made of silicon with periodically arranged pores (circular holes of radius $0.25$ carved in a periodic fashion), as depicted in~\Cref{fig:crystal_domain}. 
The difference in the domain material, silicon for the crystal and air for the pore, is encoded in the discontinuous diffusion coefficient $a = 1 / \varepsilon$, where $\varepsilon$ is the permittivity with $\varepsilon = 12.1$ in silicon (crystal) and $\varepsilon =  1$ in air (pores) \cite{Joannopolous2008}. 
We note that the units are chosen such that the speed of light is $c=1$; see \cite{Kozon2023,CorbijnSchlottbomVegtVos2024} for similar example configurations.
We set $\beta=-1$, $f=0$, and define the boundary source term as an incoming plane wave that propagates in the $x$-direction, namely $g  = \cplu \kappa(1 - (1,0)^T \cdot \hat{n}) e^{-\cplu  \kappa x}$, with $\hat{n}$ being the unit outward normal vector,
and constant wavenumber $\kappa=\omega$. Since it is of interest to model how light propagates and interacts with the pores within a photonic crystal, we take the wave frequency as a multiple of the length of a unit cell, specifically $\kappa=1$ in this section.

With no analytical solution available for this example, we show the relative $L^2$--error decay $\| \uGround - \uACMSh \|_{\LO}/\| \uGround \|_{\LO}$ of the ACMS approximation against a high-order FEM solution $\uGround$ with degree of approximation $p=10$, which we refer to as the ground truth solution.

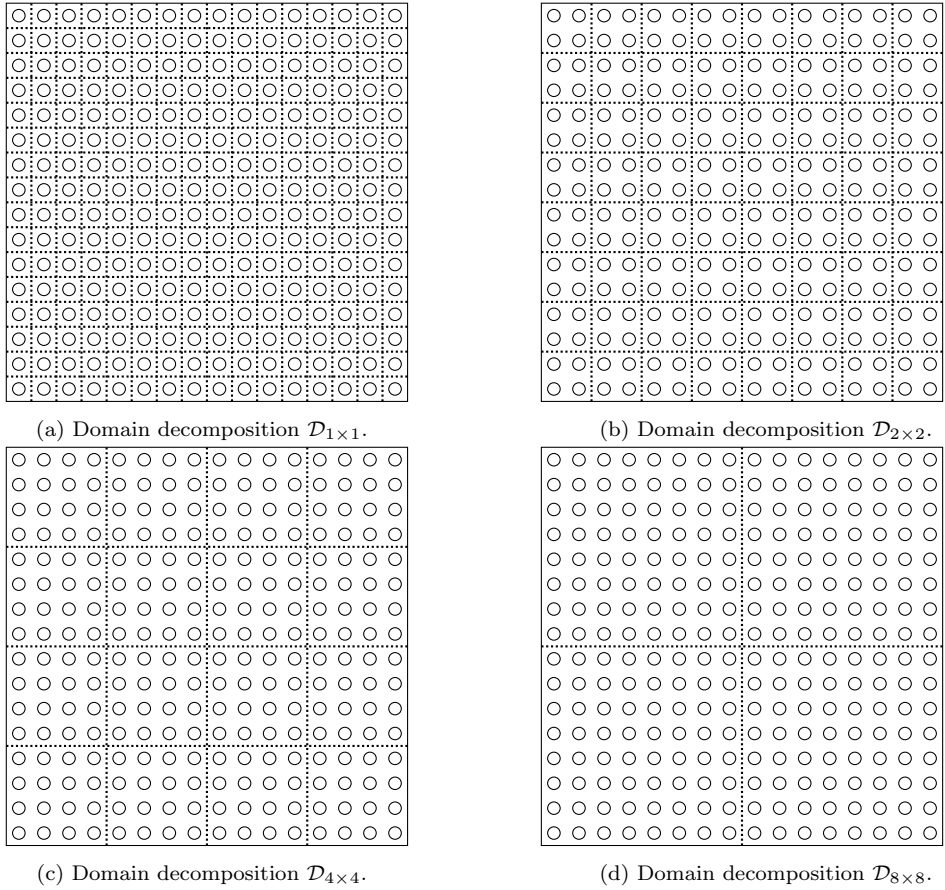
\begin{figure}[ht]
\centering 
\begin{subfigure}[b]{0.45\textwidth}
\centering
\begin{tikzpicture}[>=latex, scale = 0.33]
  \def\Ncell{16}
  \def\incl{1}
  \pgfmathsetmacro\I{0.5/\incl}
  \pgfmathsetmacro\N{\Ncell/\incl-1}
  \def\r{0.25}
  \fill[white, fill opacity=0.2, draw = black] (0,0) rectangle (\Ncell,\Ncell);
  \foreach \x in {1, ..., \N} {
  
  \begin{scope}[densely dotted, line width = 0.8]
      \draw[black] (\x * \incl, 0) -- (\x * \incl, \Ncell);
      \draw[black] (0, \x * \incl) -- (\Ncell, \x * \incl);
  \end{scope}
  }
      \foreach \i in {1, ..., \Ncell} {
        \foreach \j in {1, ..., \Ncell} {
              \draw[] (-0.5+\i,-0.5+\j) circle (\r);
        }
      }
\end{tikzpicture}
\caption{Domain decomposition $\Dincl{1}$. \label{fig:crystal_inclusions_1}}
\end{subfigure}
\hfill
\begin{subfigure}[b]{0.45\textwidth}
\begin{tikzpicture}[>=latex, scale = 0.33]
  \def\Ncell{16}
  \def\incl{2}
  \pgfmathsetmacro\I{0.5/\incl}
  \pgfmathsetmacro\N{\Ncell/\incl-1}
  \def\r{0.25}
  \fill[white, fill opacity=0.2, draw = black] (0,0) rectangle (\Ncell,\Ncell);
  \foreach \x in {1, ..., \N} {
  \begin{scope}[densely dotted, line width = 0.8]
      \draw[black] (\x * \incl, 0) -- (\x * \incl, \Ncell);
      \draw[black] (0, \x * \incl) -- (\Ncell, \x * \incl);   
  \end{scope}
  }
  \foreach \i in {1, ..., \Ncell} {
    \foreach \j in {1, ..., \Ncell} {
          \draw[] (-0.5+\i,-0.5+\j) circle (\r);
    }
  }
\end{tikzpicture}
\caption{Domain decomposition $\Dincl{2}$. \label{fig:crystal_inclusions_2}}
\end{subfigure}
\\
\centering
\begin{subfigure}[b]{0.45\textwidth}
\centering
\begin{tikzpicture}[>=latex, scale = 0.33]
  \def\Ncell{16}
  \def\incl{4}
  \pgfmathsetmacro\I{0.5/\incl}
  \pgfmathsetmacro\N{\Ncell/\incl-1}
  \def\r{0.25}
  \fill[white, fill opacity=0.2, draw = black] (0,0) rectangle (\Ncell,\Ncell);
  \foreach \x in {1, ..., \N} {
  \begin{scope}[densely dotted, line width = 0.8]
      \draw[black] (\x * \incl, 0) -- (\x * \incl, \Ncell);
      \draw[black] (0, \x * \incl) -- (\Ncell, \x * \incl);
  \end{scope}
  }
  \foreach \i in {1, ..., \Ncell} {
    \foreach \j in {1, ..., \Ncell} {
          \draw[] (-0.5+\i,-0.5+\j) circle (\r);
    }
  }
\end{tikzpicture}
\caption{Domain decomposition $\Dincl{4}$. \label{fig:crystal_inclusions_4}}
\end{subfigure}
\hfill
\begin{subfigure}[b]{0.45\textwidth}
\begin{tikzpicture}[>=latex, scale = 0.33]
  \def\Ncell{16}
  \def\incl{8}
  \pgfmathsetmacro\I{0.5/\incl}
  \pgfmathsetmacro\N{\Ncell/\incl-1}
  \def\r{0.25}
  \fill[white, fill opacity=0.2, draw = black] (0,0) rectangle (\Ncell,\Ncell);
  
  \foreach \x in {1, ..., \N} {
  \begin{scope}[densely dotted, line width = 0.8]
      \draw[black] (\x * \incl, 0) -- (\x * \incl, \Ncell);
      \draw[black] (0, \x * \incl) -- (\Ncell, \x * \incl);
  \end{scope}
  }
  \foreach \i in {1, ..., \Ncell} {
    \foreach \j in {1, ..., \Ncell} {
          \draw[] (-0.5+\i,-0.5+\j) circle (\r);
    }
  }
\end{tikzpicture}
\caption{Domain decomposition $\Dincl{8}$. \label{fig:crystal_inclusions_8}}
\end{subfigure}
\caption{$\Omega$ is a $16 \times 16$ square, decomposed into $256$ (\subref{fig:crystal_inclusions_1}), $64$ (\subref{fig:crystal_inclusions_2}), $16$ (\subref{fig:crystal_inclusions_4}) and $4$ (\subref{fig:crystal_inclusions_8}) subdomains. \label{fig:crystal_inclusions}}
\end{figure}

Denoting a $1 \times 1$ square that has a circular pore in its center a unit cell, a natural choice for a domain decomposition is to assign to each unit cell a subdomain $\Oj$ for $j=1,\ldots,256$, see \Cref{fig:crystal_inclusions_1}. 
The ACMS method is, however, applicable for more general domain decompositions, and we next investigate the accuracy and efficiency of the corresponding ACMS approximation. 
Therefore, we repeat the numerical tests for different subdivisions of the domain $\Omega$ as depicted in \Cref{fig:crystal_inclusions}, which correspond to domain decompositions into $J = 1, 4, 16, 64, 256$ subdomains with number of edges $|\E| =4, 12, 40, 144, 544$, and number of vertices $|\V| = 4, 9, 25, 81, 289$. We refer to each domain decomposition as $\Dincl{j}$, where $j^2$ indicates the number of unit cells per subdomain. As $j$ increases, the number of subdomains decreases.
The length of the edges, which influences the error estimate \eqref{eq:L2_error_bound_kappa_dep}, varies for each domain decomposition.
Hence, in order to counterbalance the effect of the increased edge length on the error, we have to consider twice as many modes in our numerical tests when going from $\Dincl{j}$ to $\Dincl{2j}$.

In~\Cref{tab:Crystal_DD_Conv}, we present the approximation error for the different domain decompositions for a fixed mesh size $h=0.05$ and for an increasing polynomial order $p=1,2,3,4,5,6$.
The number of edge modes $\Ie$, indicated in brackets, increases with the order $p$ and doubles if the subdomain size doubles.
From the data shown in the table, we observe that all domain decompositions $\Dincl{j}$, $j=1,2,4,8,16$ give similar results.
Hence, we may conclude that the constants in the error estimate \eqref{eq:L2_error_bound_kappa_dep} depend only weakly on the chosen domain decomposition. 
The slight increase in the error for $p=5$ might be explained by the fact that the respective ACMS spaces for $p=4,5,6$, which employ the same $\Ie$ if $\Dincl{j}$ is fixed, are not nested.

We next discuss the influence on the computational costs when going from $\Dincl{j}$ to $\Dincl{2j}$ and at the same time doubling the number of modes for accuracy reasons, see \eqref{eq:L2_error_bound_kappa_dep}. 
We first observe that the corresponding number of subdomains $J$ decreases by $4$ and the number of degrees of freedom on an edge $\NFEMe$ roughly doubles, and $\NFEMj=\NFEMe^2$ approximately increases by four. 
We recall \eqref{eq:acms_total_costs} 
\begin{align*} 
    \underbrace{\O(J \Ie \Cext)
    }_{\text{basis computation}}
    +
    \underbrace{\O(J \NFEMj \Ie^2 )
    }_{\text{assembly}} 
    +
    \underbrace{\O( J^2 \Ie^3)
    }_{\text{solve}}.
\end{align*}
The basis computation cost depends on the cost $\Cext$ for the extension. Assuming $\Cext=\O(J (\NFEMe^4+\Ie \NFEMe^3))$, cf. \Cref{rem:extension}, the resulting basis computation cost of $\O(J^2 (\Ie \NFEMe^4 + \Ie^2 \NFEMe^3))$ favors larger $J$, i.e., smaller $j$. 
Similarly, the assembly cost favors smaller $\Ie$, i.e., larger $J$, since $J \NFEMj$ is constant. 
The cost for solving of $\O(J^2 \Ie^3)$ favors smaller $J$, but if it scales only like $O(J^{3/2}\Ie^3)$, cf. \Cref{fig:scaling}, it becomes independent of the domain decomposition, because $J^{3/2}\Ie^3$ stays constant. 
From these considerations, we may conclude that choosing $j$ as small as possible is beneficial. 
However, parallelization allows to reduce the influence of the basis computation and assembly costs.
Moreover, if subdomains have moderate number $\NFEMj$ of degrees of freedom, then cache effects and dense linear algebra may reduce the influence of $\Cext$ on the runtime.
Since the applicability of sparse direct solvers may be limited by the size of the corresponding linear system, keeping the overall system size small is favorable.
Since the system size corresponds to the dimension of the ACMS space $\NACMS=\O(J\Ie)$, it is favorable to keep $J$ moderate, i.e., to employ larger subdomains.
For these reasons, we continue with the domain decomposition $\Dincl{2}$.
\begin{table}
\centering
\scalebox{0.65}{
\centering\small\setlength\tabcolsep{0.55em}
\begin{tabular}{ c c c c c c c}
\toprule
&&& $p$   \\
\cmidrule{2-7} 
& $1$  & $2$ & $3$ & $4$ & $5$ & $6$  \\
\toprule
$\Dincl{1}$ & $6.9 {\cdot} 10^{-2} \ (16)$ & $2.4 {\cdot} 10^{-4} \ (32)$ & $6.5 {\cdot} 10^{-6} \ (32)$ &  $4.1 {\cdot} 10^{-7} \ (64)$  & $9.1 {\cdot} 10^{-7} \ (64)$& $7.2 {\cdot} 10^{-7} \ (64)$
\\
$\Dincl{2}$ & $7.4 {\cdot} 10^{-2} \ (32)$ & $2.6 {\cdot} 10^{-4} \ (64)$ & $6.5 {\cdot} 10^{-6} \ (64)$ &  $4.1 {\cdot} 10^{-7} \ (128)$  & $8.8 {\cdot} 10^{-7} \ (128)$ & $7.0 {\cdot} 10^{-7} \ (128)$ 
\\
$\Dincl{4}$ & $7.7 {\cdot} 10^{-2} \ (64)$ & $2.8 {\cdot} 10^{-4} \ (128)$ & $6.5 {\cdot} 10^{-6} \ (128)$ &  $4.1 {\cdot} 10^{-7} \ (256)$  &  $8.7 {\cdot} 10^{-7} \ (256)$ &  $7.0 {\cdot} 10^{-7} \ (256)$ 
\\
$\Dincl{8}$ & $7.9 {\cdot} 10^{-2} \ (128)$ & $2.9 {\cdot} 10^{-4} \ (256)$ & $6.6 {\cdot} 10^{-6} \ (256)$ &  $4.1 {\cdot} 10^{-7} \ (512)$  & $8.7 {\cdot} 10^{-7} \ (512)$ & $7.0 {\cdot} 10^{-7} \ (512)$ 
\\
$\Dincl{16}$ & $8.0 {\cdot} 10^{-2} \ (256)$ & $2.9 {\cdot} 10^{-4} \ (512)$ & $6.6 {\cdot} 10^{-6} \ (512)$ &  $4.1 {\cdot} 10^{-7} \ (1024)$  & $8.6 {\cdot} 10^{-7} \ (1024)$ & $6.9 {\cdot} 10^{-7} \ (1024)$
\\
\bottomrule
\end{tabular}
}
\caption{Example in \Cref{test:numex2}, on the square domain $\Omega$ with different decompositions depicted in \Cref{fig:crystal_domain,fig:crystal_inclusions}: $\Dincl{1} $, $\Dincl{2} $, $\Dincl{4}$, $\Dincl{8} $, $ \Dincl{16} $ . The $L^2$--relative error $\| \uGround - \uACMSh \|_{\LO}/\| \uGround \|_{\LO}$ is computed against the ground truth solution. We choose mesh size $h=0.05$, degree of approximation $p=1,2,3,4,5,6$, and number of edge modes $\Ie$ as indicated in brackets.
\label{tab:Crystal_DD_Conv}
}
\end{table}

\pgfplotstableread[]{WavenumberSweepCrystalDincl2.out}\WavenumberSweepCrystal
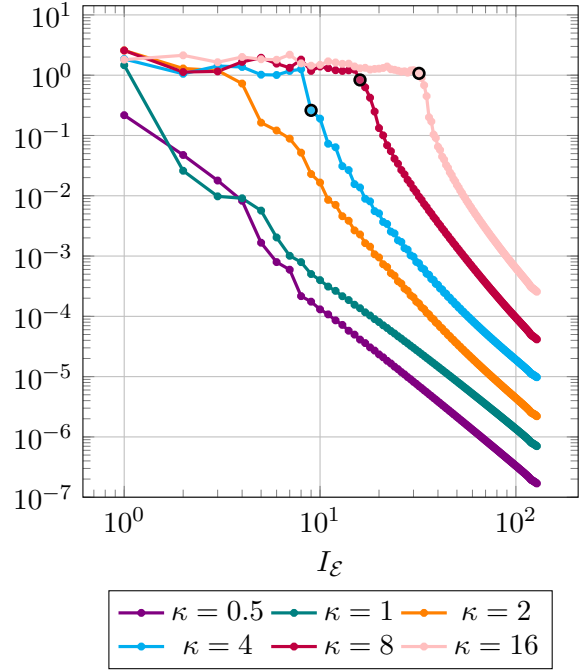
\begin{figure}[ht]
\centering\small \setlength\tabcolsep{0.55em}
\begin{tikzpicture}
  \begin{axis}[    
    width=0.6\textwidth,height=0.6\textwidth,
    ymode=log,
    xmode=log,
    legend style={
      legend columns=3,
      at={(0.5,-0.35)},
      anchor=south,
      },
    xlabel = \( \Ie \), 
    ymin = 1.0e-7, 
    xmajorgrids=true, ymajorgrids=true,
    scaled y ticks = false,
    cycle list name=philipcolors,
    ]

   \pgfplotsset{cycle list shift=15}    
    
    \addplot table[x=Ie,y expr={\thisrowno{1}} ]{\WavenumberSweepCrystal}; 
    \addplot table[x=Ie,y expr={\thisrowno{2}} ]{\WavenumberSweepCrystal};  
    \addplot table[x=Ie,y expr={\thisrowno{3}} ]{\WavenumberSweepCrystal}; 
    \addplot table[x=Ie,y expr={\thisrowno{4}} ]{\WavenumberSweepCrystal}; 
    \addplot table[x=Ie,y expr={\thisrowno{5}} ]{\WavenumberSweepCrystal}; 
    \addplot table[x=Ie,y expr={\thisrowno{6}} ]{\WavenumberSweepCrystal}; 
    
    \addplot[ color=black, mark=*, fill = white, fill opacity=0.2, mark size = 2pt , mark options={line width=1pt}]  coordinates {(9, 0.2610918804000535)};
    \addplot[  color=black, mark=*, fill = white, fill opacity=0.2, mark size = 2pt , mark options={line width=1pt}]     coordinates {(16, 0.8347304327539421)};
    \addplot[  color=black, mark=*, fill = white, fill opacity=0.2, mark size = 2pt  , mark options={line width=1pt}]  coordinates {(32, 1.0731605762685514)};


      \legend{$\kappa = 0.5$,$\kappa = 1$,$\kappa = 2$,$\kappa = 4$,$\kappa = 8$,$\kappa = 16$}
      
      \end{axis}
\end{tikzpicture}
\caption{Example in \Cref{test:numex2} on the square domain $\Omega$ with domain decomposition $\Dincl{2}$, \Cref{fig:crystal_inclusions_2}. The $L^2$--relative error of the ACMS approximation is computed against a ground truth solution, for $h= 0.05$ and $p = 6$. We show the convergence rate for an increasing number of edge modes per edge $\Ie = 1, ..., 128$ for different wavenumber values $\kappa = 0.5, 1, 2, 4, 8, 16$. 
\label{fig:Crystal_Wavenumber}}
\end{figure}
Finally, besides the dependence on the chosen decomposition, we further want to confirm the linear correlation between the number of modes and the wavenumber that was shown in~\Cref{fig:Circle_Wavenumber}, in case of a discontinuous diffusion coefficient and a different setting. In this test, we choose the wavenumbers based on an appropriate scaling of the example in \cite{CorbijnSchlottbomVegtVos2024}: wavelengths between $400$nm and $5000$nm correspond to $\omega$ values ranging from $0.628$ to $7.85$, hence we choose $\kappa \in [0.5,8]$. We additionally consider $\kappa = 16$ to test the method's limits.
In~\Cref{fig:Crystal_Wavenumber}, we show the relative $L^2$--error for the domain decomposition $\Dincl{2}$, order of approximation $p=6$ and an increasing number of edge modes $\Ie = 1, ...,128$. 
Again, we mark with a black circle the drop after which the error decreases with the predicted rate, and we observe that the number of edge modes required in order to have a monotonic error decay varies linearly with the wavenumber, i.e., when doubling $\kappa$ we have to double $\Ie$.

\subsection{Example 3: Large crystal non periodic pores} \label{test:numex3}
\begin{figure}
    \centering
      \begin{tikzpicture}[>=latex, scale = 0.5]
          \def\Nx{10}
          \pgfmathsetmacro\Nxh{\Nx/2}
          \def\Ny{10}
          \pgfmathsetmacro\Nyh{\Ny/2}
          \pgfmathsetmacro\Nya{\Ny/2-1}
          \pgfmathsetmacro\Nyb{\Ny/2+2}
          \def\r{0.25}
       
           \fill[green, fill opacity=0.2, draw = black] (0,0) rectangle (\Nx,\Ny);
          
           \draw[line width = 2, draw = cyan] (0, {\Ny/2 - 3}) --node[left] {$\gamma_{in}$}(0, {\Ny/2 + 3}) ;

           \draw[line width = 2, draw = magenta] (\Nx, {\Ny/2 - 3}) --node[right] {$\gamma_{out}$} (\Nx, {\Ny/2 + 3}) ;

           \foreach \x in {1, ..., \Nxh} {
          \draw[black!60!white, densely dotted] (\x * 2, 0) -- (\x * 2, \Ny);
      }          
      \foreach \y in {1, ..., \Nyh} {
          \draw[black!70!white, densely dotted] (0, \y * 2) -- (\Nx, \y * 2);
      }       

           \fill[pattern={north east lines}, pattern color = green!80!black] (0, \Ny/2-1) rectangle (\Nx, \Ny/2+1);
      
      \foreach \i in {1, ..., \Nx} {
        \foreach \j in {1, ..., \Nya} {
              \fill[red!30!white, opacity = 1, draw = black] (-0.5+\i,-0.5+\j) circle (\r);
        }
         \foreach \j in {\Nyb, ..., \Ny} {
              \fill[red!30!white, opacity = 1, draw = black] (-0.5+\i,-0.5+\j) circle (\r);
        }
      }
    \end{tikzpicture}
    \caption{Geometry of a crystal (green) with pores (red) and a waveguide like structure (depicted by the vertical line pattern) for $\Jx=10$ and $\Jy = 10$. Further we plot two  lines  $\gamma_{in} = [-\Jx/2] \times [-3,3]$ and $\gamma_{in} = [\Jx/2] \times [-3,3]$ in cyan and magenta, respectively. The dotted lines depict the structure of the domain decomposition used for the ACMS computations.}
    \label{fig:waveguide}  
\end{figure}
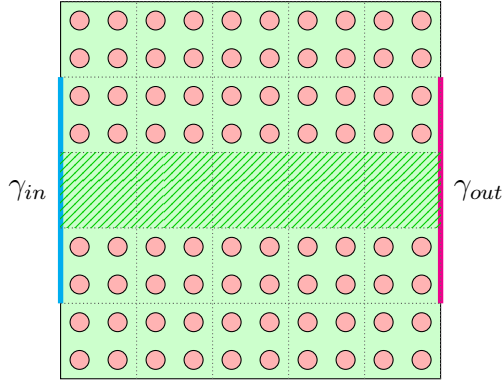

In this section we want to show that the ACMS method can be employed to simulate large but finite-size crystals, which possibly contain defects.
To do so, we compute the transmission of a "localized" wave that is placed in front of a waveguide-like crystal structure. 
We consider a rectangular shaped crystal $\Omega = [-\Jx/2, \Jx/2] \times [-\Jy/2, \Jy/2]$ consisting of $\Jx \times \Jy$ unit cells of size $ [0,1] \times [0,1]$.
We choose the same parameter as before, i.e., $f=0$, $\beta=-1$, $c=1$, $a = 1$ in the pores and $a = 1/12.1$ outside, and the right hand side
\begin{align*}
    g =  
    \begin{cases}
        \cplu \kappa \, e^{-\cplu \kappa x} \, e^{-y^2}  & \text{on } [-\Jx/2] \times [-\Jy/2, \Jy/2], \\
        0 & \text{else}.
    \end{cases}
\end{align*}
Different to the previous examples, some unit cells consist of silicon only ($a=1/12.1$), while others include pores with a radius $r = 0.25$ centered within the cell as before. 
We refer to \Cref{fig:waveguide} for an example of the geometry for $\Jx = 10$, $\Jy = 10$, where two vertical lines, $\gamma_{in} = [-\Jx/2] \times [-3,3]$ and $\gamma_{out} = [\Jx/2] \times [-3,3]$, are also added on the boundary since they are used for later computations. 
The dashed lines depict the $\Dincl{2}$ decomposition of the crystal that is considered for all computations. Furthermore, we choose polynomial order $p = 5$, mesh size $h=0.1$ and number of edge modes per edge $\Ie=16$. 

\begin{figure}
    \begin{center}
    \includegraphics[scale=1]{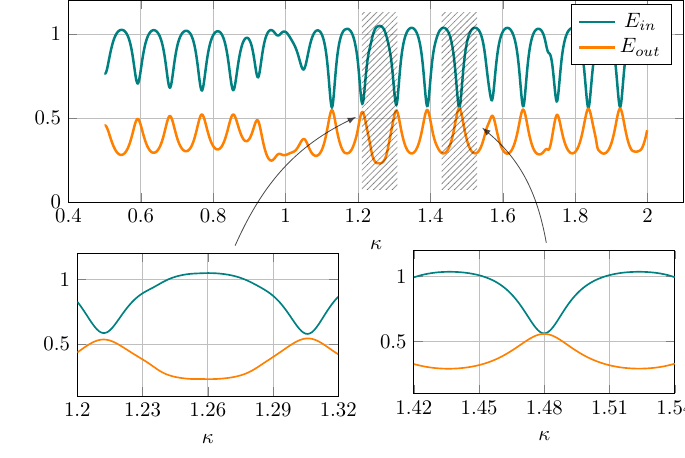}
    \end{center}
    \caption{The norms $E_{in}$ and $E_{out}$ on a crystal structure with $\Jx = \Jy = 10$ and wavenumber $\kappa \in [0.5,2]$.}
    \label{fig:sweep}
\end{figure}

\begin{figure}
    \centering
    \begin{tikzpicture}
            \node[] at (0.0,0) {\includegraphics[width = 5.2cm]{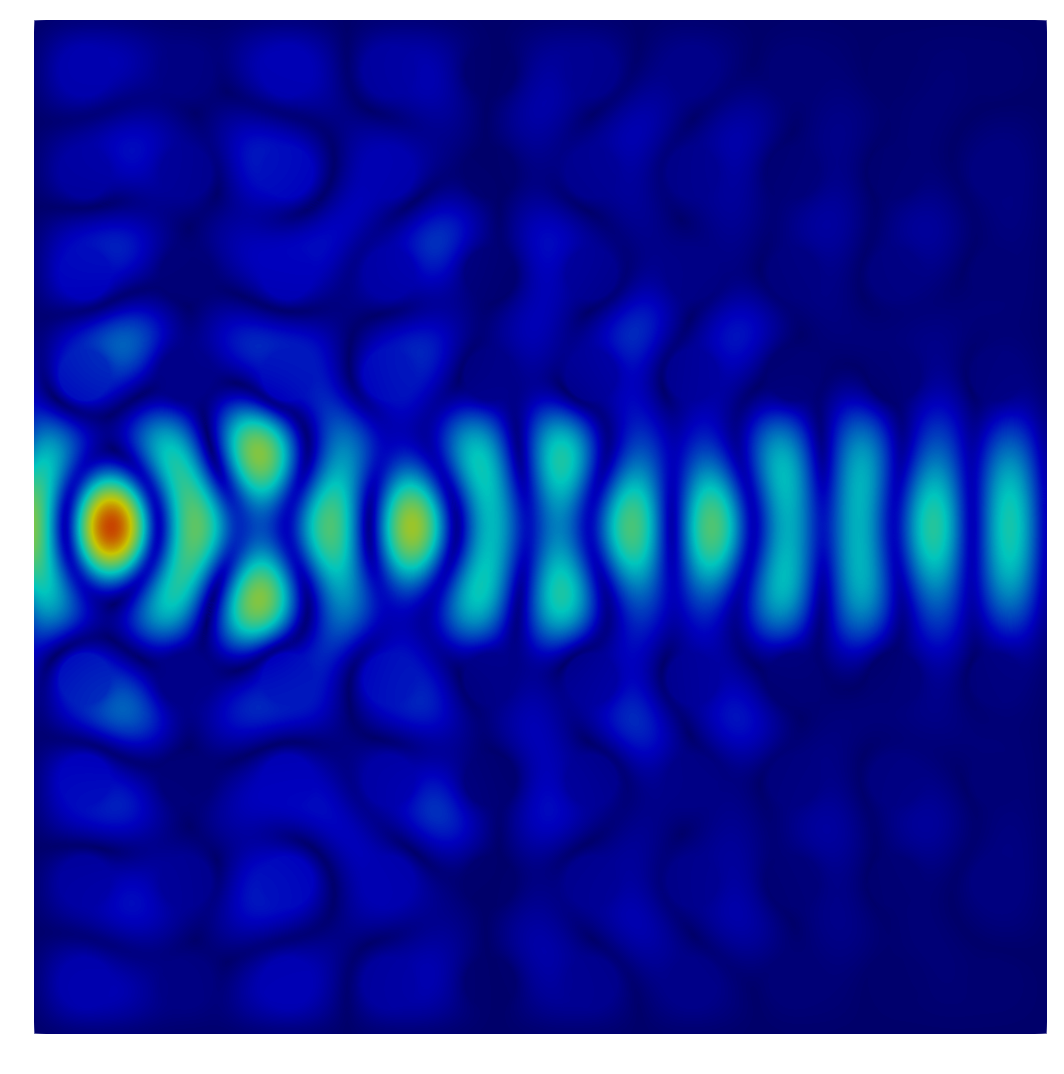}};
            \node[] at (5.5,0) {\includegraphics[width = 5.2cm]{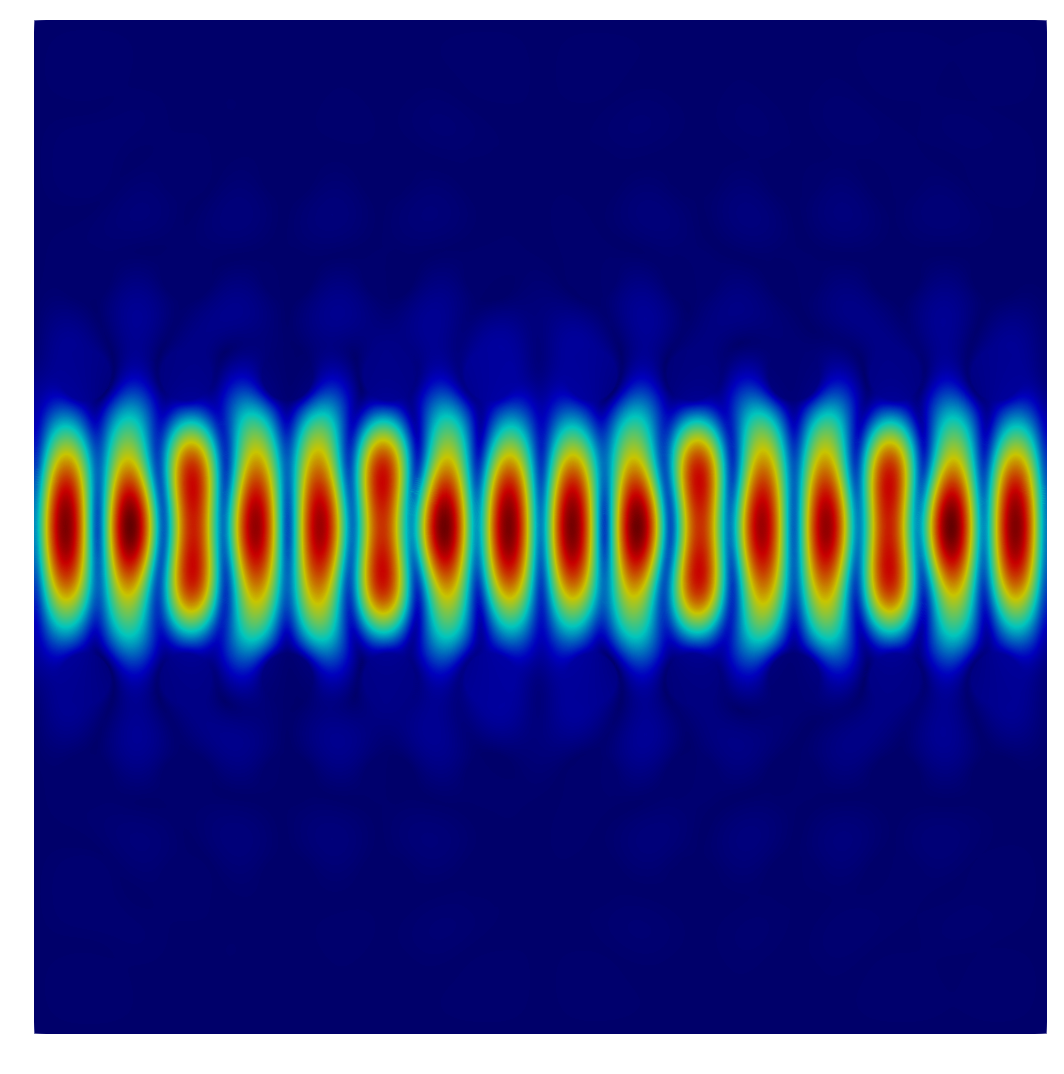}};
      \node[anchor=south west] (bar1) at (8.5,-2.8) {
    \pgfplotscolorbardrawstandalone[colormap/jet,colorbar sampled,point meta min=0,point meta max=1.8, 
                colorbar style={samples=19, width=0.2cm, height = 5.1cm, ytick={0,0.2, 0.4, 0.6, 0.8, 1.0, 1.2, 1.4, 1.6, 1.8}, 
                scaled y ticks=false, 
                yticklabel style={style={font=\footnotesize}, /pgf/number format/fixed,  /pgf/number format/precision=1}}]};

    \end{tikzpicture}
    \caption{Absolute value $|\uACMSh|$ for  $\kappa=1.26$ (left) and $\kappa = 1.48$ (right) for a crystal structure with $\Jx = \Jy = 10$.}
    \label{fig:transmission}
\end{figure}

First, we are going to perform a frequency sweep over $\kappa \in [0.5,2]$.
To compare the results we measure the $L^2$--norm on the left and right side of the crystal
\begin{align*}
    E_{in} := \int_{-\Jy/2}^{\Jy/2} | \uACMSh(-\Jx/2, y)|^2 \ dy
    \quad 
    \text{and}
    \quad 
    E_{out} := \int_{-\Jy/2}^{\Jy/2} | \uACMSh(\Jx/2, y)|^2 \ dy.
\end{align*}
Note that $E_{out}$ is related to the time-averaged energy flow through $\gamma_{out}$ via the complex Poynting vector \cite{Joannopolous2008}, while $E_{in}$ corresponds to the time-averaged energy flow through $\gamma_{in}$ up to contributions that are due to the inhomogeneous boundary terms.
In the upper plot of \Cref{fig:sweep}, where we have chosen $\Jx=\Jy=10$, one can clearly see the wavenumber dependence of $E_{in}$ and $E_{out}$ showing that the transmission pattern varies with wavenumber; see also \cite{CorbijnSchlottbomVegtVos2024} for an analogous behavior.
While some areas exhibit a transmission pattern, i.e., $E_{in}\approx E_{out}$, others suppress transmission, i.e., $E_{out}\ll E_{in}$, e.g., for $\kappa \approx 1$.
In addition, we further plot in more detail the values of $E_{in}$ and $E_{out}$ for $\kappa \in [1.2,1.32]$ and $\kappa \in [1.42, 1.54]$. First, we observe that the smallest amplitude of $E_{out}$ occurs at $\kappa = 1.26$ (within the considered wavenumber range $[0.5, 2]$), indicating poor wave transmission. In contrast, at $\kappa = 1.48$, we have $E_{in} = E_{out}$, suggesting that, at this frequency, the incoming localized wave can pass through the crystal structure, or waveguide, with minimal dissipation. Indeed, in  \Cref{fig:transmission} the absolute value of $\uACMSh$ is plotted and we can observe the expected behaviour. 

We note that the dimension of the ACMS space is $\NACMS = 516$ for $\Jx=\Jy=10$, while the underlying FEM space has a dimension $\NFEM = 253\,926$, i.e., about $\NFEMj\approx 2\,500$ degrees of freedom per unit cell. For inclusions with simple geometries and moderate wavenumbers such a resolution might not be necessary and finite element methods with lower resolution may still yield accurate solutions for the full crystal problem. 
However, high resolutions guarantee that the ACMS basis functions are realized accurately. Moreover, highly resolved FEM spaces do not cause computational issues for the ACMS method, because the local systems related to the extension operator remain small.
Next we illustrate this observation by increasing the crystal size.

\begin{figure}
    \centering
    \includegraphics[scale=1]{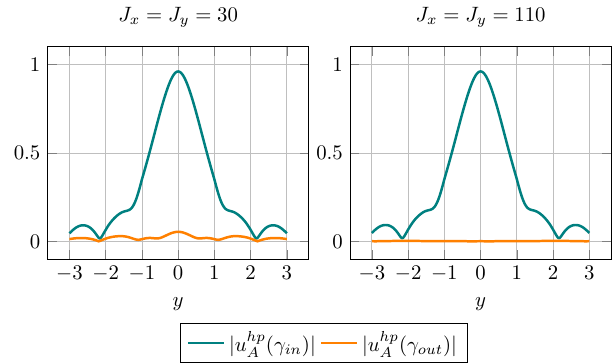}
    \caption{Absolute value $|\uACMSh|$ evaluated on $\gamma_{in}$ and $\gamma_{out}$ for a crystal with $\Jx=\Jy = 30$ (left) and $\Jx=\Jy=110$ (right) and a frequency $\kappa = 1.26$.
    \label{fig:transition_onetwosix}}
\end{figure}

\begin{figure}
    \centering
    \includegraphics[scale=1]{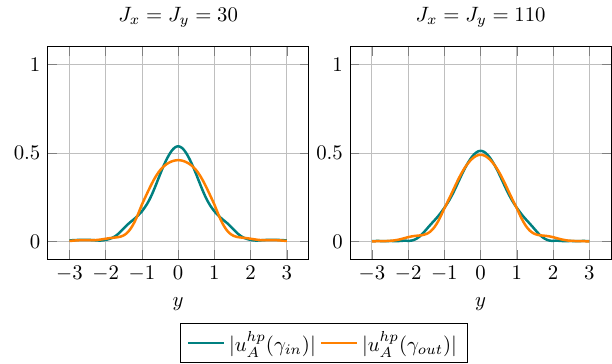}
    \caption{Absolute value $|\uACMSh|$ evaluated on $\gamma_{in}$ and $\gamma_{out}$ for a crystal with $\Jx=\Jy = 30$ (left) and $\Jx=\Jy=110$ (right) and a frequency $\kappa = 1.48$.
    \label{fig:transition_onefoureight}}
\end{figure}

Since the previous computations were performed on a very small crystal, we conclude this section by extending the computations to larger crystal structures. For this purpose, we select two frequencies, $\kappa = 1.26$ and $\kappa = 1.48$, and compute the solution on two crystals, one with $\Jx = \Jy = 30$ and the other with $\Jx= \Jy = 110$, i.e., $900$ and $12\,100$ unit cells, respectively.
We note that the dimension of the approximation spaces is $\NACMS = 4\,096$ and $\NFEM = 2\,405\,776$ for $\Jx=\Jy=30$; while $\NACMS = 52\,416$ and $\NFEM =30\,636\,426$ for $\Jx=\Jy=110$, which prohibits the solution of the FEM discretization with standard sparse direct solvers on our machines. 
In \Cref{fig:transition_onetwosix,fig:transition_onefoureight}, we plot the absolute value of the solution along the two aforementioned lines, $\gamma_{in}$ and $\gamma_{out}$. Motivated by the findings from the previous example, we expect that for $\kappa = 1.26$ the wave is less transmitted through the crystal, while for $\kappa = 1.48$ we expect that the wave is passing through the crystal with a minimal amount of absorption and possible reflection. This behavior is indeed observed in \Cref{fig:transition_onetwosix} and  \Cref{fig:transition_onefoureight} for both cases $\Jx= \Jy = 30$ and $\Jx= \Jy = 110$, where in the latter case the amplitude on $\gamma_{out}$ for $\kappa = 1.26$ has almost vanished. Additionally, it is clear that the amplitude on $\gamma_{in}$ is significantly higher for $\kappa = 1.26$ compared to $\kappa=1.48$, due to greater reflection.
The same conclusion can also be drawn for the values $E_{in}$ and $E_{out}$ in \Cref{fig:transmission} for $\kappa = 1.26$ and $\kappa = 1.48$.

\section{Discussion and Conclusions}\label{sec:discussion}

We studied the computational and implementational aspects of the approximate component mode synthesis (ACMS) method. 
The numerical study has been carried out using the Python interface of the open source software NGSolve~\cite{NGSolve}, which allows to access state-of-the-art implementations of high-order finite element methods in a non-intrusive manner.
We identified similarities of the ACMS method to usual $hp$--FEM, which permit efficient assembly and solution procedures, using sparse direct solvers. 
We gave a full computational complexity analysis, which we verified numerically. 
From this analysis, it also became clear that the system assembly can be parallelized straight-forwardly.

We study numerically the influence of the underlying $hp$--FEM approximation on the convergence of the ACMS method, focusing on varying wavenumbers. 
We find a strong improvement in accuracy when the degree of approximation is larger than one, which is in line with the discussion about the pollution effect, see \Cref{Sec:Introduction}.
We numerically confirmed the theoretical results of \cite{GiammatteoHeinleinSchlottbom2024} for different wavenumbers, see \eqref{eq:L2_error_bound_kappa_dep}, i.e., that the $L^2$--error decays cubically in the number of edge modes. 
We observed that, for increasing wavenumber, the error starts to decrease monotonically for a number of modes that scales linearly in the wavenumber. 
This is better than expected by the theoretical results given in \cite{GiammatteoHeinleinSchlottbom2024}.
Since the error decay is initially very rapid, we may conclude that the ACMS method can achieve engineering relevant accuracy by employing a moderate number of basis functions.
If higher accuracy is required, one may use exponentially convergent schemes \cite{ChenHouWang2023,MaAlberScheichl2023,Peterseim2017}. We leave a detailed computational comparison of the different methods for future research.

We further studied the dependence of the error on the employed domain decomposition, and we find a weak dependence.
Hence, we may conclude that the ACMS method offers flexibility in choosing the domain decomposition to minimize computational cost, e.g., to exploit fast dense linear algebra.
This flexibility, in turn, allowed us to simulate a large but finite-size crystal that has a line defect.

\FloatBarrier
\bibliographystyle{elsarticle-num}
\bibliography{references}

\end{document}

%% file: dofnumbering.tex
\def\Nx{2}
\def\Ny{2}
\pgfmathsetmacro\Nxm{\Nx-1}
\pgfmathsetmacro\Nym{\Ny-1}
\def\l{2}
\def\em{2}

\draw[edge] (0,0) rectangle (\Nx * \l, \Ny * \l);

\foreach \i in {1, ..., \Nx} {
\foreach \j in {1, ..., \Ny} {
  \pgfmathsetmacro\im{\i - 1};
  \pgfmathsetmacro\jm{\j - 1};
  \draw[edge] (\im * \l , \jm * \l) -- (\i * \l , \jm * \l);      
  \draw[edge] (\im * \l , \jm * \l) -- (\im * \l , \j * \l);      
}}

\foreach \i in {0, ..., \Nxm} {
\foreach \j in {0, ..., \Ny} {
  \pgfmathsetmacro\nr{int(\i + \j * \Nx)};
  \pgfmathsetmacro\dd{int(
  \nr + int(\nr/\Ny) + \nr * \em + int(\nr/\Ny) * (\Ny + 1) * \em
  )};
  \pgfmathsetmacro\dde{int(\dd + \em + 1)};
  \node[vertex] at (\i * \l , \j * \l) {\dd};
  \node[vertex] at (\i * \l + \l , \j * \l) {\dde};
}
}

\foreach \i in {0, ..., \Nxm} {
\foreach \j in {0, ..., \Ny} {
  \pgfmathsetmacro\nr{int(\i + \j * \Nx)};
  \pgfmathsetmacro\dd{int(
  \nr + int(\nr/\Ny) + \nr * \em + int(\nr/\Ny) * (\Ny + 1) * \em
  )}; 
  \foreach \k in {1,...,\em}{
    \pgfmathsetmacro\ee{int(\dd+\k)}
    \pgfmathsetmacro\offset{\l/(\em+1) * \k}
    \node[vertexedge] at (\i * \l + \offset, \j * \l) {\ee};
  }
}
}

\foreach \i in {0, ..., \Nxm} {
\foreach \j in {0, ..., \Nym} {
  \pgfmathsetmacro\nr{int(\i + \j * \Nx)};
  \pgfmathsetmacro\dd{int(
    (int(\nr/\Ny)+1) * (\Ny+1) + (int(\nr/\Ny) + 1) * \Ny * \em + int(\nr/\Ny) * (\Ny + 1) * \em + Mod(\nr,\Ny) * \em -1 
  )}; 
  \foreach \k in {1,...,\em}{
    \pgfmathsetmacro\ee{int(\dd+\k)}
    \pgfmathsetmacro\eee{int(\dd+\k+\em)}
    \pgfmathsetmacro\offset{\l/(\em+1) * \k}
    \node[vertexedge] at (\i * \l, \j * \l + \offset) {\ee};
    \node[vertexedge] at (\i * \l + \l, \j * \l + \offset) {\eee};
  }
}
}